\input ./style/arxiv-general.cfg
\documentclass[aap,MSNbibl,seceqn,dvips]{arximspdf}
\makeatletter
   \@ifpackageloaded{graphicx}{}{\usepackage{graphicx}}
\makeatother

\doi{10.1214/14-AAP1063} 
\volume{25}
\issue{5}
\pubyear{2015}
\firstpage{2867}
\lastpage{2908}
\docsubty{FLA}

\makeatletter
\newcommand{\UC}{\mathrm{UC}}
\newcommand{\rrvert}{\vert}
\newcommand{\rrVert}{\Vert}
\newcommand{\llvert}{\vert}
\newcommand{\llVert}{\Vert}
\renewcommand{\mid}{|}
\newtheorem{Theorem}{Theorem}[section]
\newproclaim{Definition}{Definition}[section]

\newproclaim{Assumption}{Assumption}[section]
\newtheorem{Lemma}{Lemma}[section]
\newproclaim{Remark}{Remark}[section]
\newproclaim{Example}{Example}[section]

\newproclaim{Definitiona}{Definition}[section]
\newtheorem{Propositiona}{Proposition}[section]
\newproclaim{Assumptiona}{Assumption}
\newtheorem{Lemmaa}{Lemma}[section]
\newproclaim{Remarka}{Remark}[section]
\newproclaim{Examplea}{Example}
\def\D{\mathbb{D}}
\def\E{\mathbb{E}}
\def\F{\mathbb{F}}
\def\H{\mathbb{H}}
\def\L{\mathbb{L}}
\def\N{\mathbb{N}}
\def\P{\mathbb{P}}
\def\R{\mathbb{R}}
\def\Ac{{\mathcal A}}
\def\Cc{{\mathcal C}}
\def\Dc{{\mathcal D}}
\def\Fc{{\mathcal F}}
\def\Hc{{\mathcal H}}
\def\Nc{{\mathcal N}}
\def\Pc{{\mathcal P}}
\def\Rc{{\mathcal R}}
\def\Tc{{\mathcal T}}
\def\Uc{{\mathcal U}}
\def\Vc{{\mathcal V}}
\def\Fh{{\widehat F}}
\def\Xh{{\widehat X}}
\def\Et{{\widetilde E}}
\def\trace{\operatorname{Tr}}
\makeatother

\begin{document}
\begin{frontmatter}

\title{Second-order BSDEs with jumps: Formulation and uniqueness\thanksref{T1}}
\runtitle{Second-order BSDEs with jumps: Formulation and uniqueness\hspace*{3pt}}

\begin{aug}
\author[A]{\fnms{Nabil}~\snm{Kazi-Tani}\ead[label=e1]{mohamed-nabil.kazi-tani@polytechnique.edu}},
\author[B]{\fnms{Dylan}~\snm{Possama\"i}\corref{}\ead[label=e2]{possamai@ceremade.dauphine.fr}\thanksref{T2}}
\and
\author[C]{\fnms{Chao}~\snm{Zhou}\ead[label=e3]{matzc@nus.edu.sg}}
\runauthor{N. Kazi-Tani, D. Possama\"i and C. Zhou}
\affiliation{Universit\'e Lyon 1, Universit\'e Paris-Dauphine\\ and National University of Singapore}
\address[A]{N. Kazi-Tani\\
ISFA\\
Universit\'e Lyon 1\\
50 Avenue Tony Garnier\\
69366 Lyon Cedex 7\\
France\\
\printead{e1}}
\address[B]{D. Possama\"i\\
Ceremade\\
Universit\'e Paris-Dauphine\\
Place du Mar\'echal de Lattre de Tassigny\\
75775 Paris Cedex 16\\
France\\
\printead{e2}}

\address[C]{C. Zhou\\
Department of Mathematics\\
National University of Singapore\\
10 Lower Kent Ridge Road\\
Block S17 Office 8--14\\
Singapore 119076\\
\printead{e3}}
\end{aug}
\thankstext{T1}{Supported in part by CMAP, Ecole Polytechnique.}
\thankstext{T2}{Supported in part by Mathematics Department at the
National University of Singapore.}

\received{\smonth{5} \syear{2013}}
\revised{\smonth{1} \syear{2014}}

%
\begin{abstract}
In this paper, we define a notion of second-order backward stochastic
differential equations with jumps
(2BSDEJs for short), which generalizes the continuous case considered
by Soner, Touzi and Zhang [\textit{Probab. Theory Related Fields}
\textbf{153} (2012) 149--190]. However, on the contrary to their formulation,
where they can define pathwise the density of quadratic variation of
the canonical process, in our setting, the compensator of the jump
measure associated to the jumps of the canonical process, which is the
counterpart of the density in the continuous case, depends on the
underlying probability measures. Then in our formulation of 2BSDEJs,
the generator of the \mbox{2BSDEJs} depends also on the underlying probability
measures through the compensator. But the solution to the 2BSDEJs can
still be defined universally. Moreover, we obtain a representation of
the $Y$ component of a solution of a 2BSDEJ as a supremum of solutions
of standard backward SDEs with jumps, which ensures the uniqueness of
the solution.
\end{abstract}

%
\end{frontmatter}

\setcounter{footnote}{2}

\section{Introduction}\label{sec1}

Motivated by duality methods and maximum principles for optimal
stochastic control, Bismut \cite{bis} studied a linear backward
stochastic differential equation (BSDE). In their seminal paper \cite
{pardpeng}, Pardoux and Peng generalized such equations to the
nonlinear Lipschitz case and proved existence and uniqueness results
in a Brownian framework. Since then, a lot of attention has been given
to BSDEs and their applications, not only in stochastic control, but
also in theoretical economics, stochastic differential games and
financial mathematics.

Given a filtered probability space $(\Omega,\mathcal F, \{
\mathcal F_t \}_{0\leq t\leq T},\mathbb P)$ generated by an
\mbox{$\mathbb R^d$-}valued Brownian motion $B$, solving a BSDE with generator
$f$ and terminal condition $\xi$ consists of finding a pair of
progressively measurable processes $(Y,Z)$ such that
%
\begin{equation}
Y_t=\xi+\int_t^T
f_s(Y_s,Z_s)\,ds-\int
_t^T Z_s \,dB_s,\qquad
\mathbb P\mbox{-a.s., }t\in[0,T]. \label{defbsde}
\end{equation}

The process $Y$ thus defined is a possible generalization of the
conditional expectation of $\xi$, since when $f$ is the null function,
we have $Y_t=\E_t^{\P} [\xi]$, and in this case, $Z$ is
the process appearing in the $(\Fc_t)$-martingale representation of
$(\E_t^{\P} [\xi])_{t\geq0}$. In the case of a filtered
probability space generated by both a Brownian motion $B$ and a Poisson
random measure $\mu$ with compensator $\nu$, the martingale
representation for $(\E_t^{\P} [\xi])_{t\geq0} $ becomes
\begin{eqnarray*}
\E^{\P}_t[\xi] &=& \mathbb E^\mathbb P[\xi]+\int_0^t Z_s \,dB_s + \int_0^t\!\int_{\R^d\setminus\{0\}}
\psi_s(x) (\mu-\nu) (dx,ds),\qquad \mathbb P\mbox{-a.s.},
\end{eqnarray*}
where $\psi$ is a predictable function.

This leads to the following natural generalization of equation (\ref
{defbsde}) to the case with jumps. We will say that $(Y,Z,U)$ is a
solution to the BSDE with jumps (BSDEJ in the sequel) with generator
$f$ and terminal condition $\xi$ if for all $t \in[0,T]$, we have
$\mathbb P$-a.s.
%
\begin{eqnarray}\label{defbsdej}
Y_t &=& \xi+\int_t^T f_s(Y_s,Z_s,U_s)\,ds-\int
_t^T Z_s \,dB_s
\nonumber\\[-8pt]\\[-8pt]\nonumber
&&{} -\int_t^T\!\!\int_{\R^d\setminus\{0\}}
U_s(x) (\mu-\nu) (dx,ds).
\end{eqnarray}

Tang and Li \cite{tangli} were the first to prove existence and
uniqueness of a solution for (\ref{defbsdej}) in the case where $f$
is Lipschitz in $(y,z,u)$. Our aim in this paper is to generalize (\ref
{defbsdej}) to the second order, as introduced recently by Soner,
Touzi and Zhang \cite{stz}. Their key idea in the definition of the
second-order BSDEs (2BSDEs) is that the equation defining the solution
has to hold $\P$-almost surely, for every $\P$ in a class of
nondominated probability measures. They then manage to prove a
uniqueness result using a representation of the solution of a 2BSDE as
an essential supremum of solutions of standard BSDEs. This
representation finds its origin in the deep link that 2BSDEs share with
stochastic control theory and PDEs. In order to shed more light on this
aspect, let us give the intuition behind this representation in the
continuous case. 

Suppose that we want to study the following fully nonlinear PDE:
%
\begin{equation}
\qquad - \frac{\partial u}{\partial t}-h \bigl(t,x,u(t,x),Du(t,x),D^2u(t,x)
\bigr)=0,
\qquad u(T,x)=g(x). \label{eq1}
\end{equation}

If the function $\gamma\mapsto h(t,x,r,p,\gamma)$ is assumed to be
convex, then it is equal to its double Fenchel--Legendre transform, and
if we denote its Fenchel--Legendre transform by $f$, we have
%
\begin{equation}
h(t,r,p,\gamma)=\sup_{a\geq0} \biggl\{\frac{1}2a\gamma
-f(t,x,r,p,a) \biggr\}. \label{eq2}
\end{equation}

From (\ref{eq2}), we expect, at least formally, that the solution $u$
of (\ref{eq1}) is going to verify
\[
u(t,x)=\sup_{a\geq0} u^a(t,x),
\]
where $u^a$ is defined as the solution of the following semi-linear PDE:
%
\begin{eqnarray}\label{eq3}
- \frac{\partial u^a}{\partial t}-\frac{1}2aD^2u^a(t,x)+f
\bigl(t,x,u^a(t,x),Du^a(t,x),a \bigr)&=&0,
\nonumber\\[-8pt]\\[-8pt]\nonumber
u^a(T,x) &=& g(x).
\end{eqnarray}

Since $u^a$ is linked to a classical BSDE, the 2BSDE associated to $u$
should correspond (in some sense) to the supremum of the family of
BSDEs indexed by $a$. Furthermore, changing the process $a$ can be
achieved by changing the probability measure under which the BSDE is
written. We also emphasize that the generator of the BSDEs depends
explicitly on $a$, which is actually the density of the quadratic
variation of the martingale driving the BSDE.

For the sake of clarity, we will now briefly outline the main
differences and difficulties due on the one hand to second-order
framework and on the other hand to our jump setting.
\begin{longlist}[(iii)]
\item[(i)] We remind the reader that our aim is to introduce an equation
similar to~(\ref{defbsdej}). But as shown above in the continuous
case, the generator $f$ will have to depend on the density of
$[B,B]^c$, the pathwise continuous part of the quadratic variation
$[B,B]$ of the canonical process $B$, but since we are in a jump
setting, it will also have to depend on the compensator of the random
jump measure $\mu_B$ associated to $B$. 
Exactly as in the continuous case, we can always give a pathwise
definition of the density of $[B,B]^c$, which gives us directly an
aggregator. However, it is generally impossible to find an aggregator
for the compensator of the jump measure; see Section~\ref{secaggregation} for more details. This forces us to consider in our
jump setting 2BSDEs whose generator depends explicitly on the
underlying probability measure. This is an important difference with
the framework considered in \cite{stz}. However, in spite of this, the
solution to the $2$BSDEJs is still-defined independently of the
probability measures considered.


\item[(ii)] A second major difference with (\ref{defbsdej}) in the
second-order case, is, as we recalled earlier, that the BSDE has to
hold $\P$-almost surely for every probability measure $\P$ lying in a
wide family of probability measures. Under each $\P$, $[B,B]^c$ and
$\mu_B$ have, respectively, a prescribed density and a prescribed jump
measure compensator. This is why we can intuitively understand the
2BSDEJ (\ref{2bsdej}) as a BSDEJ with model uncertainty, where the
uncertainty affects both the quadratic variation and the jump measure
of the process driving the equation. 

\item[(iii)] The last major difference with (\ref{defbsdej}) in the
second-order case is the presence of an additional nondecreasing
process $K$ in the equation. To have an intuition for~$K$, one has to
have in mind representation (\ref{representationref}) that we prove in
Theorem~\ref{uniqueref}, stating that the $Y$ part of a solution of a
2BSDEJ is an essential supremum of solutions of standard BSDEJs. The
process $K$ maintains $Y$ above any solution $y^{\P}$ of a BSDEJ, with
given quadratic variation and jump measure under $\P$. The process $K$
is then formally analogous to the nondecreasing process appearing in
reflected BSDEs (as defined in \cite{elkarkap}, e.g.).
\end{longlist}

There are many other possible approaches in the literature to handle
volatility and/or jump measure uncertainty in stochastic models \cite
{alp,denis,dhp,bionkerv}. Among them, Peng
\cite{peng1} introduced a notion of Brownian motion with uncertain
variance structure, called $G$-Brownian motion. This process is defined
without making reference to a given probability measure. It refers
instead to the $G$-Gaussian law, defined by a partial differential
equation; see \cite{peng2} for a detailed exposition and references.
We also would like to mention the very recent works by Neufeld and Nutz
\cite{nn2,nn3}, which appeared during the revision of this paper,
which provide an elegant and very important extension to the work of
Peng, by allowing a very general type of uncertainties for the whole
triplet of characteristics of a given L\'evy process, very in much in
the spirit of the approach we follow in this paper. We emphasize that
their approach is very general, which is why they have to deal with
delicate measurability issues, and that, roughly speaking, their
results could be used to define the solution to a 2BSDEJ with a
generator equal to $0$, without having to impose any continuity
assumptions on the terminal condition.

Finally, recall that Pardoux and Peng \cite{pardpeng} proved that if
the randomness in $f$ and $\xi$ is induced by the current value of a
state process defined by a forward stochastic differential equation,
the solution to the BSDE (\ref{defbsde}) could be linked to the
solution of a semilinear PDE by means of a generalized Feynman--Kac
formula. Similarly, Soner, Touzi and Zhang \cite{stz} showed that
2BSDEs generalized the point of view of Pardoux and Peng, in the sense
that they are connected to the larger class of fully nonlinear PDEs.
In this context, the 2BSDEJs are the natural candidates for a
probabilistic solution of fully nonlinear integro-differential
equations. This is the purpose of our accompanying paper \cite{kpz}.

The rest of this paper is organized as follows. In Section~\ref
{section1}, we introduce the set of probability measures on the
Skorohod space $\D$ that we will work with. Using the notion of
martingale problems on $\D$, we construct probability measures under
which the canonical process has given characteristics. 
In Section~\ref{sec2BSDE}, we define the notion of 2BSDEJs and show
how it is linked to standard BSDEJs. Section~\ref{section2} is
devoted to our uniqueness result and some a priori estimates. The
\hyperref[app]{Appendix} is dedicated to the proof of some important technical results
needed throughout the paper.

\section{Preliminaries} \label{section1}

\subsection{A primer on 2BSDEJs and main difficulties}
Before giving all notation in detail and a precise definition of
$2$BSDEJs, we would like to start by presenting the main object of
interest in this paper, as well as the main difficulties we need to
address in our framework.

First, as mentioned in the \hyperref[sec1]{Introduction}, we shall consider the
following $2$BSDEJ, for $0\leq t\leq T$ and $\P$-a.s.:
\begin{eqnarray*}
Y_t &=& \xi+\int_t^T
\widehat{F}^{\mathbb P}_s(Y_s,Z_s,U_s)\,ds
-\int_t^T Z_s\,dB^{\mathbb P,c}_s
\\[-1pt]
&&{} -
\int_t^T\!\!\int_{E}
U_s(x) \tilde{\mu}^{\mathbb P}_B(dx,ds) +
K^{\mathbb P}_T-K^{\mathbb P}_t, 
\end{eqnarray*}
for every $\P\in\mathcal{P}^\kappa_H$, which is a family, not
necessarily dominated, of local martingale probability measures. These
different probability measures represent the model uncertainty. $B^{\P
,c}$ and $ \tilde{\mu}^{\mathbb P}_B$ denote, respectively, the
continuous local martingale part and the compensated jump measure
associated to the purely discontinuous local martingale part of the
canonical process $B$ under any local martingale measure $\P$. We
reiterate that in contrast to (\ref{defbsdej}), we have to add a
nondecreasing process $K^\P$ to account for the fact that solutions
to 2BSDEJs have to be understood as suprema of families of classical BSDEJs.

Let us now highlight the new difficulties in our framework compared to
the continuous $2$BSDEs as considered in \cite{stz}. While a crucial
issue in their definition of the 2BSDEs is the aggregation of the
quadratic variation of the canonical process $B$ under a wide family of
probability measures, here, in general, the aggregation of the jump
compensators associated to $B$ is not\vspace*{1pt} possible; see Section~\ref
{secaggregation} for more details. This is the reason why the
generator $\widehat{F}^\P$ and the compensated jump measure
$\widetilde\mu^\P$ above depend explicitly on the probability
measure, through the jump compensator defined under each $\P$. This is
an important difference which may lead one to think that it might not
be possible to define the solution $(Y,Z,U)$ of a 2BSDEJ in a universal
way (i.e., to say so that it does not depend explicitly on the
measure~$\P$). This would be very unfortunate from the point of view
of applications, since, if we look, for instance, at classical problems
of portfolio optimization in finance, the process $Z$ is usually
related to the corresponding optimal investment strategy. Therefore, in
a context of uncertainty, one will definitely need an optimal strategy
which works for every possible model, that is to say for every measure
$\P$. Nonetheless, we prove that the solution of a $2$BSDEJ,
$(Y,Z,U)$, can still be constructed in such a way that it is defined
for all $\omega$, independently of probability measures; we refer the
reader to our companion paper \cite{kpz} for more details.

Another crucial point in the definition of 2BSDEs in \cite{stz}, is
that they work under a set of measure corresponding to the so-called
strong formulation of stochastic control. Roughly speaking, this
corresponds to considering the laws under the Wiener measure of
stochastic integrals with respect to the canonical process $B$, with
the constraint that these integrands have to take values in the space
of symmetric definite positive matrices. Such a choice has several
extremely important advantages: first of all, it allows them to define
their measures through a unique reference measure (i.e., the Wiener
measure), and even more importantly, they showed that all the measures
thus constructed satisfy the martingale representation property and the
Blumenthal 0--1 law, which are known to be fundamental properties for
the wellposedness of classical BSDEs (which, as recalled in the
\hyperref[sec1]{Introduction} are a kind of nonlinear martingales). Therefore, in our
framework, we have to be able to retrieve the strong formulation.
However, there is no longer any clear choice for a reference measure as
soon as jumps are added into the mix. This will therefore lead us to
consider a whole family of reference measures, which makes in turn the
problem more complicated. Moreover, even though the only assumption
needed in \cite{stz} to retrieve the martingale representation
property is that the admissible volatilities are symmetric definite
positive (and therefore invertible), in a framework with jumps, we need
to consider special jumps compensators with some restrictions, but
which still are flexible enough to be able to model as many types of
jump measure uncertainty as possible. Overcoming these main
difficulties is the most important contribution of this paper.

\subsection{The stochastic basis}
We first introduce the notations used in the paper. Let $\Omega:= \D
([0,T],\mathbb R^d)$ be the space of c\`adl\`ag paths defined on
$[0,T]$ with values in $\R^d$ and such that $w(0)=0$, equipped with
the Skorohod topology, so that it is a complete, separable metric
space; see \cite{bil}, for instance.

We denote $B$ the canonical process, $\mathbb F:= \{\mathcal
F_t \}_{0\leq t\leq T}$ the filtration generated by~$B$, $\mathbb
F^+:= \{\mathcal F_t^+ \}_{0\leq t\leq T}$ the right limit
of $\mathbb F$ and for any $\mathbb P$, $\overline{\mathcal
F}_t^\mathbb P:=\mathcal F_t^+\vee\mathcal N^\mathbb P(\mathcal
F_t^+)$ where
\[
\mathcal N^\mathbb P(\mathcal G):= \bigl\{E\in\Omega,\mbox{ there exists
$\widetilde E\in\mathcal G$ such that $E\subset\widetilde E$ and
$\mathbb P(\widetilde E)=0$} \bigr\}.
\]

We then define as in \cite{stz} a local martingale measure $\mathbb P$
as a probability measure such that $B$ is a $\mathbb P$-local martingale.
We then associate to the jumps of $B$ a counting measure $\mu_{B}$,
which is a random measure on $\mathbb R^+\times E$ equipped with its
Borel \mbox{$\sigma$-}field $\mathcal B(\R^+)\times\mathcal B(E)$ (where
$E:=\mathbb R^d\setminus\{0\}$), defined pathwise by
%
\begin{equation}
\mu_{B} \bigl(A,[0,t] \bigr):= \sum_{0<s\leq t}
\mathbf{1}_{\{\Delta B_s \in A\}
}\qquad\forall t \geq0, \forall A \in\mathcal B(E).
\end{equation}

We recall that (see, e.g., Theorem I.4.18 in \cite{jac}) under
any local martingale measure $\P$, we can decompose $B$ uniquely into
the sum of a continuous local martingale, denoted by $B^{\P,c}$, and a
purely\vspace*{1pt} discontinuous local martingale, denoted by $B^{\P,d}$. We
emphasize that such a decomposition depends on the underlying
probability measure. Then we define $\overline{\mathcal P}_W$ as the
set of all local martingale measures $\mathbb P$, such that $\mathbb P$-a.s.:
\begin{longlist}[(ii)]
\item[(i)] The quadratic variation of $B^{\P,c}$ is
absolutely continuous with respect to the Lebesgue measure $dt$, and
its density takes values in $\mathbb S^{>0}_d$, which is the space of
all $d\times d$ real valued positive definite matrices.
\item[(ii)] The compensator $\lambda^\mathbb P_t(dx,dt)$ of
the jump measure $\mu_B$ exists under $\mathbb P$ and can be
decomposed, for some $\F$-predictable random measure $\nu^\P$ on
$E$, as follows:
\[
\lambda^\P_t(dx,dt)=\nu^\P_t(dx)\,dt.
\]
\end{longlist}

We will denote by $\widetilde\mu_{B}^\mathbb P(dx,dt)$ the
corresponding compensated measure, and for simplicity, we will often
call $\nu^\P$ the compensator of the jump measure associated to $B$.

%
\begin{Remark}
In this paper, we always assume that under the probability measures
that we consider, the canonical process is a local martingale, whose
quadratic variation and jump compensator change depending on the
measure considered. Formally, it means that we do not consider drift
uncertainty. Hence, the reader may wonder why we do not consider more
generally probability measures under which the canonical process is a
semimartingale with a triplet of characteristics which can all vary.
In a nutshell, the framework considered here is completely sufficient
for us in order to give wellposedness results for 2BSDEs with jumps,
and we did not want to make our presentation confusing. However, we
emphasize that all the above results can be easily extended to the more
general case of drift, volatility and jump uncertainty. For related
results, the reader may consult \cite{nutz2}, and the very recent
preprint \cite{nn2}.
\end{Remark}

In this discontinuous setting, we will say that a probability measure
$\mathbb P\in\overline{\mathcal P}_W$ satisfies the martingale~representation property if for any $(\overline{\mathbb F}^\mathbb
P,\mathbb P)$-local martingale $M$, there exists a unique $\overline{\mathbb F}^\mathbb P$-predictable processes $H$ and a unique
$\overline{\mathbb F}^\mathbb P$-predictable function $U$ such that
$(H,U)\in\mathbb H^2_{\mathrm{loc}}(\mathbb P)\times\mathbb J^2_{\mathrm{loc}}(\mathbb
P)$ (these spaces are defined later in Section~\ref{secspace}) and
\[
M_t=M_0+\int_0^tH_s\,dB_s^{\P,c}+
\int_0^t\!\int_EU_s(x)
\widetilde\mu_{B}^\mathbb P(dx,ds),\qquad \mathbb P\mbox{-a.s.}
\]

We now follow \cite{stz3} and introduce the so-called universal
filtration. For this we let $\mathcal P$ be a given subset of
$\overline{\mathcal P}_W$ and define the following:

\begin{Definition}
A property is said to hold $\mathcal P$-quasi-surely ($\mathcal P$-q.s.
for short), if it holds $\mathbb P$-a.s. for all $\mathbb P\in\mathcal P$.
\end{Definition}

%

\subsection{Aggregation (or the absence of it)}\label{secaggregation}
In this section, we discuss issues related to aggregation of the
quadratic variation of the canonical process $B$ and the absence of
aggregation of the jump compensators associated to the jumps of $B$.\vspace*{1pt}

Let $\Pc\subset\overline{\mathcal P}_W$ be a set of nonnecessarily
dominated probability measures, and let $\{X^{\P}, \P\in\Pc\}$
be a family of random variables indexed by $\Pc$. One can think, for
example, of the stochastic integrals $X_t^{\P}:= {^{(\mathbb P)}}\!\int
_0^t H_s \,dB_s$, where $\{H_t, t \geq0 \}$ is a predictable process.


\begin{Definition}
An \textit{aggregator} of the family $\{X^{\P}, \P\in\Pc\}$ is
a r.v. $\Xh$ such that
\begin{eqnarray*}
\Xh&=& X^{\P},\qquad\P\mbox{-a.s., for every } \P\in\Pc.
\end{eqnarray*}
\end{Definition}

Bichteler \cite{bich}, Karandikar \cite{kar}, and more recently Nutz
\cite{nutz} all showed in different contexts and under different
assumptions,\vspace*{1pt} that it is possible to find an aggregator for the It\^o
stochastic integrals ${^{(\P)}}\! \int_0^t H_s \,dB_s$. A direct
consequence of this result is the possibility to define the quadratic
variation process $\{[B,B]_t, t\geq0\}$ pathwise.\footnote{The
following construction was proposed to us by Marcel Nutz, whom we thank
warmly. It is also used in a more general context in the recent
preprint \cite{nn2}, where the absence of aggregation in a jump
setting is also made clear. We urge the reader to consult their very
interesting results.} Indeed, using It\^o's formula, we can write for
any local martingale measure~$\P$,
\begin{eqnarray*}
[B,B]_t &=& \llvert B_t\rrvert^2 - 2 \int
_0^t B_{s^-}\,dB_s,\qquad \P\mbox{-a.s.},
\end{eqnarray*}
and the aggregation of the stochastic integrals automatically yields
the aggregation of the bracket $\{[B,B]_t, t\geq0\}$.

Next, since $[B,B]$ has finite variation, we can define its
path-by-path continuous part $[B,B]^c$ (by subtracting the sum of the
jumps) and finally the corresponding density
\[
\hat a_t:=\mathop{\overline{\lim}}_{\varepsilon\downarrow0} \frac{[B,B]^c_t-
[B,B]^c_{t-\varepsilon}}{\varepsilon}.
\]
%
Notice that since for any local martingale measure $\P$,
\[
[B,B]^c= \bigl\langle B^{\P,c} \bigr\rangle,\qquad \P\mbox{-a.s.},
\]
then $\hat a$ coincides with the density of quadratic variation of
$B^{\P,c}$, $\P$-a.s. Therefore $\hat a$ takes values in $\mathbb
S^{>0}_d$, $dt\times d\P$-a.e., and
\[
\hat a_t= \frac{d\langle B^{\P,c}\rangle_t}{dt},\qquad \P\mbox{-a.s.}
\]

More generally than the above examples, Soner, Touzi and Zhang \cite
{stz3}, motivated by the study of stochastic target problems under
volatility uncertainty, obtained an aggregation result for a family of
probability measures corresponding to the laws of some continuous
martingales on the canonical space $\Omega= \Cc(\R^+,\R^d)$, under
a \textit{separability} assumption on the quadratic variations (see
their Definition~4.8) and under an additional \textit{consistency}
condition (which is usually only necessary) for the family to
aggregate. A related result, not limited to the case of volatility
uncertainty was then obtained by Cohen \cite{cohen}. In our setting,
this naturally leads to the question of whether it is possible or not
to find an aggregator for the family of jump compensators $\nu^\P$.

However, unlike with the quadratic variation which can be either
obtained through the Doob--Meyer decomposition of the local
submartingale $\langle B \rangle$ or through the use of It\^o's
formula, the predictable compensator can only be obtained thanks to the
Doob--Meyer decomposition of the nondecreasing process $[B,B]$. It is
therefore obvious that this compensator depends explicitly on the
underlying probability measure, and it is not clear at all whether an
aggregator always exists or not.

This actually goes deeper, and in any reasonable setting of jump
uncertainty, it is actually not possible to define such an aggregator,
as showed in the following simple examples.

%
\begin{Example}
Consider two probability measures $\P_1$ and $\P_2$ such that under
$\P_1$ the canonical process $B$ is a L\'evy process with
characteristics $(0,1,\lambda_1\delta_{\{1\}})$ where the intensity
of jumps $\lambda_1$ is a constant, and under $\P_2$ the canonical
process $B$ is a L\'evy process with characteristics $(0,1,\lambda
_2\delta_{\{1\}})$ where $\lambda_2$ is a constant different from
$\lambda_1$ (it is a classical result that these probabilities are
uniquely defined on the Skorohod space $\D$). Since only the jump
intensities are different, from the classical theory of change of
measures, we know that $\P_1$ and $\P_2$ are actually equivalent, so
their null-sets are the same, and we cannot find an aggregator which
is simultaneously equal to $\lambda_1$ and $\lambda_2$ on the same
support of $\P_1$ and $\P_2$. 
\end{Example}

%
\begin{Example}
Even in the case of pure jump martingale measures, we can still have a
simple counterexample. Actually, we consider probability measures such
that the canonical process $B$ is a pure jump local martingale. Under
$\P_1$, $B$ is a L\'evy process with characteristics $(0,0,2\delta_{\{
1\}}+4\delta_{\{-1\}})$, and under $\P_2$, $B$ is a L\'evy process
with characteristics $(0,0,3\delta_{\{1\}}+5\delta_{\{-1\}})$. Then
$\P_1$ and $\P_2$ are equivalent, and they are both martingale measures.
\end{Example}

Therefore, we will not try to aggregate the family of compensators of
jump measure in our formulation of 2BSDEJs. We emphasize that this
feature is also shared by the drift, which can obviously be changed by
using Girsanov's theorem. Hence, among the three elements of the
characteristic triplet of a semimartingale, as defined in \cite{jac},
for instance, the quadratic variation plays a peculiar role, in the
sense that this is the only one which can be aggregated when
uncertainty about this triplet is considered.

Notwithstanding this unavoidable fact, as proved in the following
sections, the solution to the 2BSDEJs, which is our object of interest,
can still be aggregated. To begin, we will use in the following
subsection the notion of martingale problem for semimartingales with
general characteristics (as defined in the book by Jacod and Shiryaev
\cite{jac} to which we refer), in order to construct a probability
measure under which the canonical process has a given quadratic
variation and a given jump measure.

\subsection{Characterization by martingale problems}
In this section, we extend the connection between diffusion processes
and probability measures established in \cite{stz3} thanks to weak
solutions of SDEs, to our general jump case with the more general
notion of martingale problems.

Let $\mathcal{N}$ be the set of $\mathbb F$-predictable random
measures $\nu$ on $\mathcal{B}(E)$ satisfying
%
\begin{eqnarray}\label{hypnu}
\int_0^t\!\int_{E}
\bigl(1 \wedge\llvert x\rrvert^2 \bigr)\nu_s(\omega,dx)\,ds &<& + \infty
\quad\mbox{and}
\nonumber\\[-8pt]\\[-8pt]\nonumber
\int_0^t\!\int_{ \llvert x\rrvert>1 } x \nu_s(\omega,dx)\,ds &<& +\infty
\qquad\mbox{for all } \omega\in\Omega,
\end{eqnarray}
and let $\Dc$ be the set of $\mathbb F$-predictable processes $\alpha
$ taking values in $\mathbb S_d^{>0}$ with
\[
\int_0^T \bigl\llvert\alpha_t(
\omega) \bigr\rrvert \,dt<+\infty\qquad\mbox{for all }\omega\in\Omega.
\]

We define a martingale problem as follows:

%
\begin{Definition}
For $\mathbb F$-stopping times $\tau_1\leq\tau_2$, for $(\alpha,\nu
)\in\mathcal D\times\mathcal N$ and for a probability measure
$\mathbb P_1$ on $\mathcal F_{\tau_1}$, we say that $\P$ is a
solution of the \textit{martingale problem} $(\P_1,\tau_1,\tau
_2,\alpha,\nu)$ if:
\begin{longlist}[(ii)]
\item[(i)] $\P= \P_1$ on $\Fc_{\tau_1}$.
\item[(ii)] The canonical process $B$ on $[\tau_1,\tau_2]$
is a semimartingale under $\P$ with characteristics
\[
\biggl( -\int_{\tau_1}^{\cdot} \int_E
x \mathbf{1}_{\llvert
x\rrvert>1}\nu_s(dx)\,ds, \int_{\tau_1}^{\cdot}
\alpha_s \,ds, \nu_s(dx)\,ds \biggr).
\]
\end{longlist}
\end{Definition}

%
\begin{Remark}\label{remmartpbm}
We refer to Theorem II.2.21 in \cite{jac} for the fact that $\P$ is
a solution of the martingale problem $(\P_1,\tau_1,\tau_2,\alpha,\nu)$
if and only if the following properties hold:
\begin{longlist}[(ii)]
\item[(i)] $\P= \P_1$ on $\Fc_{\tau_1}$.

\item[(ii)] The processes $M$, $J$ and $Q$ defined below are $\P
$-local martingales on $[\tau_1,\tau_2]$:
\begin{eqnarray}
M_t &:=& B_t - \sum_{\tau_1 \leq s \leq t}
\mathbf{1}_{\llvert
\Delta B_s\rrvert>1} \Delta B_s + \int_{\tau_1}^t
\int_E x \mathbf{1}_{\llvert x\rrvert>1}\nu_s(dx)\,ds,\qquad
\tau_1 \leq t \leq\tau_2,\nonumber
\\
J_t &:=& M_t^2 - \int_{\tau_1}^{t}
\alpha_s \,ds - \int_{\tau_1}^{t} \int
_E x^2 \nu_s(dx)\,ds,\qquad
\tau_1 \leq t \leq\tau_2,\nonumber
\\
Q_t &:=& \int_{\tau_1}^{t} \int
_E g(x)\mu_B(dx,ds) - \int
_{\tau
_1}^{t} \int_E g(x)
\nu_s(dx)\,ds,\nonumber
\\
\eqntext{\tau_1 \leq t \leq\tau_2,\forall g \in\Cc^+ \bigl(\R^r \bigr),}
\end{eqnarray}
where $\Cc^+(\R^r)$ is a discriminating family of bounded Borel
functions; see Remark~II.2.20 in \cite{jac} for more details.
\end{longlist}
\end{Remark}

We say that the martingale problem associated to $(\alpha,\nu)$ has a
unique solution if, for every stopping times $\tau_1, \tau_2$ and for
every probability measure $\P_1$, the martingale problem $(\P_1,\tau
_1,\tau_2,\alpha,\nu)$ has a unique solution.

Let now $\overline{\Ac}_W$ be the set of $(\alpha, \nu) \in\Dc
\times\Nc$, such that there exists a solution to the martingale
problem $(\P^0,0,+\infty,\alpha,\nu)$, where $\P^0$ is such that
$\P^0(B_0=0)=1$.

We also denote by $\Ac_W$ the set of $(\alpha, \nu) \in\overline
{\Ac}_W$ such that there exists a unique solution to the martingale
problem $(\P^0,0,+\infty,\alpha,\nu)$. We denote $\P^{\alpha
}_{\nu}$ this unique solution and finally set
\[
\Pc_W:= \bigl\{ \P^{\alpha}_{\nu}, (\alpha, \nu)
\in\Ac_W \bigr\}.
\]

%
\begin{Remark}
We take here as an initial condition that $B_0=0$. There is actually no
loss of generality, since at the end of the day, the probability
measures under which we are going to work will all satisfy the
Blumenthal 0--1 law. Hence, $B_0$ will have to be a constant, and we
choose $0$ for simplicity.
\end{Remark}

\subsection{The strong formulation}
We now face the following problem. As reminded in the \hyperref[sec1]{Introduction}, the
predictable martingale representation property is a crucial ingredient
for the whole BSDE theory, as well as the Blumenthal 0--1 law. Hence
the set $\Pc_W$ defined above is far too large for our purpose. In
this section, we will concentrate on a subset of $\mathcal P_W$ which
will only contain probability measures that do satisfy the predictable
representation property and the Blumenthal 0--1 law. For this purpose,
let us first consider any so-called L\'evy measure, that is to say any
deterministic (i.e., which does not depend on $\omega$) measure $F\in
\mathcal N$. It is a well-known result that for any such measure $F$,
$(I_d,F)\in\mathcal A_W$, and that the corresponding unique solution
$\mathbb P_{0,F}:=\mathbb P_F^{I_d}$ satisfies the predictable
martingale representation property as well as the Blumenthal 0--1 law.
Let us then define
\[
\mathcal A_{\mathrm{det}}:= \bigl\{(I_d,F), F\in\mathcal N\mbox{ and $F$ is deterministic} \bigr\}.
\]
We would like to use this set as a base to build a class of probability
measures under which the canonical process has, formally, the following
dynamics:
%
\begin{equation}
\label{alphabeta}dB_t=\alpha_s^{1/2}\,dW_s+
\int_E\beta_s(x) \bigl(\mu(dx,ds)-F(dx)\,ds
\bigr),
\end{equation}
for some given processes $\alpha$ and $\beta$, and where $W$ is a
Brownian motion and $\mu$ is a Poisson random measure with compensator
$F$. This can usually be done by considering the law under $\P_{0,F}$
of a well chosen stochastic process; see (\ref{Xalphabeta}) below.
There are then two questions one should ask:
\begin{longlist}[(ii)]
\item[(i)] How large can one choose the corresponding family of
compensators $F$ while ensuring that the predictable martingale
representation property and the Blumenthal 0--1 law hold?
\item[(ii)] Since, on a fundamental level, the notion of a 2BSDEJ that
we want to define corresponds to a stochastic control problem where
the objective function is a family of BSDEJs indexed by the family of
probability measures considered, the chosen class of controls (i.e.,
here the family of compensators $F$) has to be rich enough for the
dynamic programming property to hold. In particular, the family of
compensators has to be stable by concatenation and bifurcation; see,
for instance,~\cite{bt}.
\end{longlist}

Since the family $\mathcal A_{\mathrm{det}}$ is clearly not stable by
concatenation and bifurcation (recall that the compensators in
$\mathcal A_{\mathrm{det}}$ are deterministic), it has to be enlarged, but
in such a way that we do not lose either the predictable martingale
representation property or the Blumenthal 0--1 law. Such a result can be
achieved by a classical construction, detailed in Section~\ref
{secdef}, by considering the so-called separable class of coefficients
generated by $\Ac_{\mathrm{det}}$ (see Definition \ref{sep}), which we
denote $\widetilde{\Ac}_{\mathrm{det}}$. We also designate by $\Pc
_{\widetilde{\Ac}_{\mathrm{det}}}$ the set of measures corresponding to
this separable class of coefficients. Then, in virtue of Proposition
\ref{toutrestevrai}, all the measures in $\Pc_{\widetilde{\Ac
}_{\mathrm{det}}}$ do\vspace*{1pt} satisfy the predictable martingale representation
property and the Blumenthal 0--1 law.

For simplicity, we let $\mathcal V$ designate the measure $F\in
\mathcal N$ such that $(I_d,F)\in\widetilde{\mathcal A}_{\mathrm{det}}$.
Moreover, we will still denote $\mathbb P_{0,F}:=\mathbb P_F^{I_d}$,
for any $F\in\mathcal V$.

Let us now detail what kind of processes $\alpha$ and $\beta$ we can
choose in (\ref{alphabeta}). For $\alpha$, we can basically take any
process in $\Dc$. For $\beta$ however, the situation is a bit more
complicated, since the admissible $\beta$ will necessarily have to
depend on the measure $F\in\Vc$ chosen. First of all, we introduce
the following set $\mathcal R_F$ of $\F$-predictable functions $\beta
\dvtx E\longmapsto\R$ such that for Lebesgue almost every $s\in[0,T]$,
%
\[
\llvert\beta_s\rrvert(\omega,x)\leq C \bigl(1\wedge\llvert x
\rrvert\bigr),\qquad F_s(dx)\mbox{-a.e., for every }\omega\in\Omega,
\]
and for every $\omega\in\Omega$,
\begin{eqnarray*}
x&\longmapsto&\beta_s(\omega,x)\mbox{ is strictly monotone on the
support of the law of }\Delta B_s\\
&&\mbox{ under }\P_{0,F}.
\end{eqnarray*}

We will then denote by $\beta^{(-1)}_s(\omega,\cdot)$ the
corresponding inverse function. While the first condition is common,
since it implies in particular that for every $\beta\in\mathcal R_F$,
we have
\[
\int_0^T\!\!\int_E\llvert
\beta_s\rrvert^2(x)F_s(dx)\,ds<+\infty,\qquad \mathbb P_{0,F}\mbox{-a.s.},
\]
the second one may seem surprising. Nonetheless, it is a natural
condition in our context, since, as we will see, it will be linked to
problems of stochastic flow inversion for SDEs with jumps; see below
for more details.

Next, for each $F\in\mathcal V$ and for each $(\alpha,\beta)\in
\mathcal D\times\mathcal R_F$, we define
\[
\mathbb P^{\alpha,\beta}_F:=\mathbb P_{0,F}\circ
\bigl(X^{\alpha,\beta}_{\cdot} \bigr)^{-1},
\]
where
%
\begin{eqnarray}\label{Xalphabeta}
X^{\alpha,\beta}_t &:=& \int_0^t
\alpha_s^{1/2}\,dB_s^{\P_{0,F},c}
\nonumber\\[-8pt]\\[-8pt]\nonumber
&&{} +\int_0^t\!\int_E
\beta_s(x) \bigl(\mu_B(dx,ds)-F_s(dx)\,ds
\bigr),\qquad\mathbb P_{0,F}\mbox{-a.s.}
\end{eqnarray}

We then define
\[
\overline{\mathcal P}_S:=\bigcup_{F\in\mathcal V}
\bigl\{\mathbb P^{\alpha,\beta}_F, (\alpha,\beta)\in\mathcal D
\times\mathcal R_F \bigr\}.
\]

%
\begin{Remark}
Let us discuss a bit the kind of measures that are in the set
$\overline{\Pc}_S$. First of all, there are almost no restrictions
(except mild integrability conditions) on the admissible $\alpha$.
This means that basically, our framework covers all types of volatility
uncertainty. However, when it comes to the jump compensators which are
allowed, the situation is more complicated. Indeed, according to a
result of Jacod (see \cite{jacod}, Theorem 14.53, page~471), if we
take one measure $F\in\Ac_{\mathrm{det}}$ which is nonatomic and with
infinite mass, then every $\nu\in\Nc$ can be written as the image of
$F$ by some $\F$-predictable function $\beta$. Therefore, it would
appear that there was no need for us to consider more than one $F$.
However, the strong formulation that we consider is tailor-made so that
we can recover the predictable martingale representation property, and
as we will see below, this puts restrictions on the possible $\beta$
we can consider (namely they have to be invertible). Hence
considering only one $F$ could seriously limit the range of
compensators we can reach. This is the reason why we chose to consider
a whole family of measures $F$. However, it is a difficult problem to
know how large the set of compensators we consider is when compared to
$\Nc$. Nonetheless, from the point of view of applications, we think
that it does not induce any important restrictions, since the set $\Vc
$ by itself contains already more than all the possible compensators of
additive processes.
\end{Remark}

Notice then that $\alpha$ is the density of the quadratic variation of
the continuous part of $X^{\alpha,\beta}$ and
\[
dB_s^{\P_{0,F},c}=\alpha_s^{-1/2}\,dX_s^{\alpha,c},
\]
under $\mathbb P_{0,F}$. Moreover, the compensator of the measure
associated to the jumps of $X^{\alpha,\beta}$ is the measure $\nu
^{F,\beta}(dx)\,dt$ where
\[
\nu^{F,\beta}_t(\omega,A):=\int_E{
\mathbf1}_{\beta_t(\omega,x)\in
A}F_t(\omega,dx)\qquad\mbox{for any }A\in\mathcal
B(E),
\]
that is to say the image of the measure $F$ by $\beta$. Besides, we
have $\Delta X_s^{\alpha,\beta}=\beta_s(\Delta B_s)$ under $\mathbb P_{0,F}$.

Before\vspace*{1pt} pursuing, we would like to be able to define for any $\P\in
\mathcal P_W$ and any $F\in\mathcal V$ a process $L^{\P,F}$, whose
law under $\P$ is the same as the law of $B$ under $\P_{0,F}$. If $F$
were deterministic, then this would amount to constructing an additive
process which would be the sum of $\P$-Brownian motion and a pure jump
$\P$-martingale with compensator $F$, which is classical result. When
$F\in\mathcal V$, it is indeed random, but it has the special
structure (\ref{doublesum}). Hence, the previous construction can
easily be carried out recursively in this case. If in addition, the
probability measure $\P=:\P^{\alpha,\beta}_F$ is actually in
$\overline{\Pc}_S$, then we can instead define
%
\begin{eqnarray}\label{Lpab}
L^{\P^{\alpha,\beta}_F,F}_\cdot:=W^{\P^{\alpha,\beta}_F}_t+
\int_0^\cdot\int_E
\beta^{(-1)}_s(x) \bigl(\mu_B(dx,ds)-
\nu^{\P^{\alpha,\beta}_F}_s (dx )\,ds \bigr),
\nonumber\\[-8pt]\\[-8pt]
\eqntext{\P^{\alpha,\beta}_F\mbox{-a.s.},}
\end{eqnarray}
where $W^{\P^{\alpha,\beta}_F}$ is a $\P^{\alpha,\beta}_F$
Brownian motion defined by
\[
W^{\P^{\alpha,\beta}_F}_t:=\int_0^\cdot\hat
a_s^{-1/2}\,dB_s^{\P
^{\alpha,\beta}_F,c},
\]
and where we remind the reader that since the law of $B$ under $\P
^{\alpha,\beta}_F$ is the same as the law of $X^{\alpha,\beta}$
under $\P_{0,F}$, the support of the law of the jumps of $L^{\P
^{\alpha,\beta}_F,F}$ under $\P^{\alpha,\beta}_F$ is the image by
$\beta$ of the support of the law of the jumps of $B$ under $\P
_{0,F}$, so that $\beta^{(-1)}$ is indeed well-defined in the above expression.

Then, $\overline{\mathcal P}_S$ is a subset of $\mathcal P_W$, and we
have by definition for any $F\in\mathcal V$,
%
\begin{eqnarray}\label{eqdistrib}
&& \mbox{the }\mathbb P^{\alpha,\beta}_F\mbox{-distribution of }
\bigl(B,\hat a, \nu^{\mathbb P^{\alpha,\beta}_F},L^{\mathbb
P^{\alpha,\beta}_F,F} \bigr)
\nonumber\\[-8pt]\\[-8pt]\nonumber
&&\qquad = \mbox{the }\mathbb
P_{0,F}\mbox{-distribution of } \bigl(X^{\alpha,\beta
},\alpha,
\nu^{F,\beta},B \bigr).
\end{eqnarray}
Let us note immediately that the above implies that $B$ has the
following characteristics under $\P^{\alpha,\beta}_F$:
%
\begin{eqnarray}\label{eqcharB}
\hat a_t(B_\cdot) &=& \alpha_t
\bigl(L^{\P^{\alpha,\beta
}_F,F}(B_{\cdot}) \bigr)\quad\mbox{and}\quad
\nonumber\\[-8pt]\\[-8pt]\nonumber
\nu_t^{\mathbb P^{\alpha,\beta
}_F}(B_\cdot,dx)&=&
\nu_t^{F,\beta} \bigl(L^{\P^{\alpha,\beta
}_F,F}(B_\cdot),dx
\bigr),\qquad \P^{\alpha,\beta}_F\mbox{-a.s.}
\end{eqnarray}

Now we want to recover the predictable martingale representation
property. One possible solution would be to have a characterization of
$\overline{\Pc}_S$ in terms of completed filtrations, exactly as in
Lemma~8.1 of \cite{stz3} in the continuous case. Roughly speaking,
their result uses crucially a fact, which translated in our notation, reads
\[
\overline{\F^B}^{\P_{0,F}}\subset\overline{\F^{X^{\alpha,\beta
}}}^{\P_{0,F}}.
\]
When there are no jump terms, this result is actually trivial, as soon
as the matrix $\alpha$ is invertible (notice that also in our case,
the reverse inclusion is immediate). However in our setting, because
the jumps of $X^{\alpha,\beta}$ and $B$ are related by
\[
\Delta X_s^{\alpha,\beta}=\beta_s(\Delta
B_s),\qquad \P_{0,F}\mbox{-a.s.},
\]
even though we know that $X^{\alpha,\beta}$ and $B$ jump at the same
times, if the function $\beta$ is not invertible on the support of the
law of the jumps of $B$ under $\P_{0,F}$, we cannot identify the size
of a jump of $B$ by only observing a jump of $X^{\alpha,\beta}$. This
is a well-known problem in the literature of SDEs in a jump setting;
see, for example, Fujiwara and Kunita \cite{fk} or Protter \cite
{protter}. This is exactly the reason why we assumed that the
invertibility of the maps $\beta\in\mathcal R_F$.


We then have the following characterization of $\overline{\Pc}_S$,
which is similar to Lemma~8.1 in \cite{stz3}:

%
\begin{Lemma}
$\overline{\mathcal P}_S= \{\mathbb P\in\mathcal P_W, \exists
F\in\mathcal V, \overline{\mathbb F^{L^{\mathbb P,F}}}^\mathbb
P=\overline{\mathbb F^B}^\mathbb P \}$.
\end{Lemma}

\begin{pf}
First of all, let $\P^{\alpha,\beta}_F\in\overline{\mathcal P}_S$.
Then by definition, we have
\[
\overline{\F^B}^{\mathbb P_{0,F}}=\overline{\F^{X^{\alpha,\beta
}}}^{\mathbb P_{0,F}}.
\]
Now\vspace*{1pt} we can use (\ref{eqdistrib}) to obtain that $\overline{\mathbb
F^{L^{\mathbb P^{\alpha,\beta}_F,F}}}^{\mathbb P^{\alpha,\beta}_F}=
\overline{\mathbb F^B}^{\mathbb P_F^{\alpha,\beta}}$.\vadjust{\goodbreak}

Conversely, let $\mathbb P\in\mathcal P_W$ be such that there exists
some $F\in\mathcal V$ and $\overline{ \F^{L^{\mathbb P,F}}}^{\mathbb
P}= \overline{\mathbb F^B}^{\mathbb P}$. Then, there exists some
measurable function $\zeta$ such that $B_{\cdot}=\zeta(L_{\cdot}^{\mathbb P,F}),
\P$-a.s.

Now notice that by definition, the law of $\zeta(B_\cdot)$ under $\P
_{0,F}$ is the same as the law of $\zeta(L_{\cdot}^{\mathbb P,F})$ under $\P
$; that is, this is the same as the law of $B$ under $\P$. Therefore,
since $B$ is a $(\mathbb P,\F)$-local martingale by definition, $\zeta
(B_\cdot)$ is a $(\P_{0,F},\F)$-local martingale. However, since, as
recalled above, $\P_{0,F}$ has the predictable martingale
representation property, there exist a $\overline{\F}^{\P
_{0,F}}$-predictable process $\alpha$ and a \mbox{$\overline{\F}^{\P
_{0,F}}$-}predictable function $\beta$ such that
\begin{eqnarray}
\zeta(B)_t=\int_0^t
\alpha_s^{1/2}\,dB_s^{\P_{0,F},c}+\int
_0^t\!\int_E
\beta_s(x) \bigl(\mu_B(dx,ds)-F_s(dx)\,ds
\bigr),\nonumber
\\
\eqntext{\mathbb P_{0,F}\mbox{-a.s.}}
\end{eqnarray}
Notice also that we can always take a $\P_{0,F}$ version of $\alpha$
and $\beta$ which is \mbox{$\F$-}predictable. Then, we actually have $\zeta
(B_\cdot)=X^{\alpha,\beta}_F$. Fix now any measurable and bounded
function $\varphi$, and we have
\begin{eqnarray*}
\mathbb E^\mathbb P \bigl[\varphi(B_\cdot) \bigr] &=& \mathbb
E^{\mathbb P} \bigl[\varphi\bigl(\zeta\bigl(L_\cdot^{\mathbb P,F}
\bigr) \bigr) \bigr]
= \mathbb E^{\mathbb P_{0,F}} \bigl[\varphi\bigl(
\zeta(B_\cdot) \bigr) \bigr]
\\
&=& \mathbb E^{\mathbb P_{0,F}} \bigl[\varphi
\bigl(X^{\alpha,\beta}_\cdot\bigr) \bigr]
= \mathbb E^{\mathbb P^{\alpha
,\beta}_F}\bigl[\varphi(B_\cdot) \bigr],
\end{eqnarray*}
which means that $\P=\P^{\alpha,\beta}_F$.
\end{pf}

As an immediate consequence, we deduce the following, since for any
$F\in\mathcal V$, we have the martingale representation property for
any $(\overline{\F^{L^{\P,F}}}^\P,\P)$-local martingale, and the
Blumenthal 0--1 law holds for the filtration $\overline{\F^{L^{\P
,F}}}^\P$.

\begin{Lemma}
Every probability measure in $\overline{\Pc}_S$ satisfies the
predictable martingale representation property and the Blumenthal 0--1 law.
\end{Lemma}

\begin{pf}
Fix some $\P^{\alpha,\beta}_F\in\overline{\Pc}_S$. Let us start
with the predictable martingale representation property. We start by
denoting for simplicity,
\[
\overline{\F}^{\alpha,\beta}:=\overline{\mathbb F^{L^{\mathbb
P^{\alpha,\beta}_F,F}}}^{\mathbb P^{\alpha,\beta}_F}.
\]
Let $M$ be a $(\overline{\F}^{\P^{\alpha,\beta}_F},\P^{\alpha,\beta
}_F)$-local martingale; then it is also a $(\overline{\F
}^{\alpha,\beta},\P^{\alpha,\beta}_F)$-local martingale. Then by
the standard predictable martingale representation theorem, we know
that there exist a unique pair $(\widetilde H,\widetilde U)$ of
$\overline{\F
}^{\alpha,\beta}$-predictable process and function such that,
$\mathbb P^{\alpha,\beta}_F$-a.s.
\begin{eqnarray*}
M_t&=& M_0+\int_0^t
\widetilde H_s\,dW_s^{\P^{\alpha,\beta}_F}
\\
&&{} +\int
_0^t\!\int_E\widetilde
U_s(x) \bigl(\mu_{L^{\P^{\alpha,\beta
}_F},F}(dx,ds)-F_s
\bigl(L^{\P^{\alpha,\beta}_F,F}(B_\cdot),dx \bigr)\,ds \bigr).
\end{eqnarray*}
Define
\[
H:=\hat a^{-1/2}\widetilde H\quad\mbox{and}\quad U(x):=\widetilde U
\bigl(\beta^{(-1)}(x) \bigr).
\]
Then, using (\ref{eqcharB}) and (\ref{Lpab}), we obtain directly that
\begin{eqnarray}
M_t=M_0+\int_0^tH_s\,dB_s^{\P^{\alpha,\beta}_F,c}+
\int_0^t\!\int_EU_s(x)
\bigl(\mu_{B}(dx,ds)-\nu^{\P^{\alpha,\beta
}_F}_s(dx)\,ds \bigr),\nonumber
\\
\eqntext{\P^{\alpha,\beta}_F\mbox{-a.s.}}
\end{eqnarray}

The Blumenthal 0--1 law can then be directly deduced; see the proof of
Lemma~8.2 in \cite{stz3} for details.
\end{pf}

\section{Preliminaries on 2BSDEJs}\label{sec2BSDE}
\subsection{The nonlinear generator}
In this subsection we will introduce the function which will serve as
the generator of our 2BSDEJ. Let us define the following spaces for
$p\geq1$:
\[
\widehat{L}^p:= \bigl\{\xi, \mathcal F_T
\mbox{-measurable, s.t. }\xi\in L^p(\nu),\mbox{ for every }\nu\in\Nc
\bigr\}.
\]

We then consider a map
\[
H_t(\omega,y,z,u,\gamma,\widetilde{v})\dvtx[0,T]\times\Omega\times
\mathbb{R}\times\mathbb{R}^d\times\widehat{L}^2 \times
D_1\times D_2\rightarrow\mathbb{R},
\]
where $D_1 \subset\mathbb{R}^{d\times d}$ is a given subset
containing $0$ and $D_2 \subset\widehat{L}^1$ is the domain of $H$ in the
variable $\widetilde{v}$.

Define the following conjugate of $H$ with respect to $\gamma$ and
$\widetilde{v}$ by
%
\begin{eqnarray}\label{Ffenchel}
\quad && F_t(\omega,y,z,u,a,\nu)
\nonumber\\[-8pt]\\[-8pt]\nonumber
&&\qquad :=\sup_{\{\gamma,\widetilde{v}\} \in D_1
\times
D_2} \biggl\{
\frac{1}2\trace(a\gamma)+\int_{E} \widetilde{v}(e)
\nu(de)-H_t (\omega,y,z,u,\gamma,\widetilde{v} ) \biggr\},
\end{eqnarray}
for $a \in\mathbb S_d^{>0}$ and $\nu\in\mathcal{N}$.

In the remainder of this paper, we formulate the needed hypothesis for
the generator directly on the function $F$, and the BSDEs we consider
also include the case where $F$ does not take the form (\ref
{Ffenchel}). Nonetheless, this particular form allows us to retrieve
easily the framework of the standard BSDEs or of the $G$-stochastic
analysis on the one hand (see Sections~\ref{connectionbsdes} and \ref
{connectionglevy}), and to establish the link with the associated PDEs
on the other hand. In the latter cases, $H$ is evaluated at $\tilde
{v}(\cdot) = Av(\cdot)$, where $A$ is the following nonlocal
operator, defined for any $\mathcal{C}^1$ function $v$ on $\R^d$ with
bounded gradient, and $y \in\R^d$ by:
\begin{eqnarray}
(Av) (t,y) (e):= v(t,e+y) - v(t,y)- \mathbf{1}_{\{\llvert e\rrvert\leq
1\}} e. (\nabla v)
(t,y)\nonumber
\\
\eqntext{\mbox{for } e \in E\mbox{ and } t \in[0,T].}
\end{eqnarray}

The assumptions on $v$ ensure that $(Av)(t,y)(\cdot)$ is an element of
$\widehat{L}^1$.

The operator $A$ applied to $v$ will not appear again in the paper, but
this particular nonlocal form comes from the intuition that the
2BSDEJs is an essential supremum of standard BSDEJs. Indeed, solutions
to Markovian BSDEJs provide viscosity solutions to some parabolic
partial integro-differential equations with similar nonlocal
operators; see \cite{bbp} for more details.

We define for any $\P\in\overline{\Pc}_S$
%
\begin{equation}
\widehat{F}_t^{\mathbb P}(y,z,u):=F_t \bigl(y,z,u,
\hat{a}_t,\nu_t^{\mathbb P} \bigr)\quad\mbox{and}\quad
\widehat{F}_t^{\mathbb
P,0}:=\widehat{F}_t^{\mathbb P}(0,0,0).
\end{equation}

We\vspace*{1.5pt} denote by $D^1_{F_t(y,z,u)}$ the domain of $F$ in $a$ and by
$D^2_{F_t(y,z,u)}$ the domain of $F$ in $\nu$, for a fixed $(t,\omega
,y,z,u)$. As in \cite{stz} we fix a constant $\kappa\in(1,2]$ and
restrict the probability measures in $\mathcal{P}_H^\kappa\subset
\overline{\mathcal{P}}_{S}$.

%
\begin{Definition}\label{def}
$\mathcal{P}_H^\kappa$ consists of all $\mathbb P \in\overline{\Pc
}_S$ such that:
\begin{longlist}[(ii)]
\item[(i)] $\mathbb E^\mathbb P [\int_0^T\!\!\int_E \llvert
x\rrvert^2{\nu}^{\mathbb
P}_t(dx)\,dt ]<+\infty$;
\item[(ii)] $\underline{a}^\mathbb P \leq\hat{a}\leq
\overline{a}^\mathbb P$, $dt\times d\mathbb P\mbox{-a.s. for some }
\underline{a}^\mathbb P, \overline{a}^\mathbb P \in\mathbb{S}_d^{>0}$
and
\[
\mathbb{E}^{\mathbb{P}} \biggl[ \biggl(\int_0^T\bigl|
\widehat{F}_t^{\mathbb P,0}\bigr|^\kappa \,dt \biggr)^{{2}/\kappa}
\biggr]<+\infty.
\]
\end{longlist}
\end{Definition}

%
\begin{Remark}
The above conditions assumed on the probability measures in $\mathcal
P^\kappa_H$ ensure that under any $\mathbb P\in\mathcal P^\kappa_H$,
the canonical process $B$ is actually a true c\`adl\`ag martingale.
This will be important when we will define standard BSDEJs under each
of these probability measures.
\end{Remark}

We now state our assumptions on the function $F$ which will be our main
interest in the sequel.

\begin{Assumption} \label{assumphref}
(i) The domains $D^1_{F_t(y,z,u)}=D^1_{F_t}$ and
$D^2_{F_t(y,z,u)}=D^2_{F_t}$ are independent of $(\omega,y,z,u)$.

(ii) For fixed $(y,z,a,\nu)$, $F$ is $\mathbb{F}$-progressively
measurable in $D^1_{F_t} \times D^2_{F_t} $.

(iii) The following uniform Lipschitz-type property holds. For
all $(y,y',z,\break z',u, t,a,\nu,\omega)$,
\[
\bigl\llvert F_t(\omega,y,z,u,a,\nu)- F_t \bigl(
\omega,y',z',u,a,\nu\bigr) \bigr\rrvert\leq C \bigl(
\bigl\llvert y-y' \bigr\rrvert+ \bigl\llvert a^{1/2}
\bigl(z-z' \bigr) \bigr\rrvert\bigr).
\]

(iv) For all $(t,\omega,y,z,u^1,u^{2},a,\nu)$, there exist two
processes $\gamma$ and $\gamma'$ such that
\begin{eqnarray*}
\int_{E}\delta^{1,2} u(x)\gamma'_t(x)
\nu(dx)&\leq& F_t \bigl(\omega,y,z,u^1,a,\nu\bigr)-
F_t \bigl(\omega,y,z,u^2,a,\nu\bigr)
\\
&\leq&\int
_{E}\delta^{1,2} u(x)\gamma_t(x)
\nu(dx),
\end{eqnarray*}
where $\delta^{1,2} u:=u^1-u^2$ and $c_1(1\wedge\llvert x\rrvert) \leq
\gamma
_t(x) \leq c_2(1\wedge\llvert x\rrvert)$ with $-1+\delta
\leq c_1\leq0,
c_2\geq0$,
and $c_1'(1\wedge\llvert x\rrvert) \leq\gamma'_t(x) \leq
c_2'(1\wedge\llvert x\rrvert)$ with $-1+\delta\leq c_1'\leq
0, c_2'\geq0$,
for some $\delta>0$.

(v) $F$ is uniformly continuous in $\omega$ for the Skorohod
topology, that is to say that there exists some modulus of continuity
$\rho$ such that for all $(t,\omega,\omega',y,z,u,a,\nu)$,
\[
\bigl\llvert F_t(\omega,y,z,u,a,\nu)-F_t \bigl(
\omega',y,z,u,a,\nu\bigr) \bigr\rrvert\leq\rho
\bigl(d_S \bigl(\omega_{\cdot\wedge t},\omega'_{\cdot\wedge t}
\bigr) \bigr),
\]
where $d_S$ is the Skorohod metric and where $\omega_{\cdot\wedge
t}(s):=\omega(s\wedge t)$.
\end{Assumption}

%
\begin{Remark}
Assumptions (i) and (ii) are classic in the second-order framework; see
\cite{stz}. Lipschitz assumption (iii) is standard in the BSDE theory
due to the paper \cite{pardpeng}. Hypothesis (iv) allows us to have a
comparison theorem in the framework with jumps; it was introduced in
\cite{roy} and is also present in \cite{bbp} in the form of an
equality. The last hypothesis (v) of uniform continuity in $\omega$ is
also proper to the second-order framework; it is linked to our
intensive use of regular conditional probability distributions in \cite
{kpz} to construct our solutions in a pathwise manner, thus avoiding
complex issues related to negligible sets. Moreover, we emphasize that
unlike \cite{stz}, we consider here the Skorohod topology instead of
the topology induced by the uniform norm. This is linked to the fact
that we need our space $\Omega$ to be separable. Furthermore, notice
that if we restrict ourselves to the Wiener space as in \cite{stz}, we
recover their assumption since the topologies induced by the uniform
norm and the Skorohod metric are then equivalent. Nonetheless, this
property will only be useful for us in our accompanying paper \cite{kpz}.
\end{Remark}

%
\begin{Remark}
(i) For $\kappa_1<\kappa_2$, applying H\"older's inequality
gives us
\[
\mathbb{E}^{\mathbb{P}} \biggl[ \biggl(\int_0^T
\bigl\llvert\widehat{F}_t^{\mathbb P,0} \bigr\rrvert
^{\kappa_1} \,dt \biggr)^{{2}/{\kappa
_1}} \biggr] \leq C\mathbb{E}^{\mathbb{P}}
\biggl[ \biggl( \int_0^T \bigl\llvert
\widehat{F}_t^{\mathbb P,0} \bigr\rrvert^{\kappa_2} \,dt
\biggr)^{{2}/{\kappa
_2}} \biggr],
\]
where $C$ is a constant. Then it is clear that $\mathcal{P}^\kappa_H$
is decreasing in $\kappa$.

(ii) Assumption \ref{assumphref}, together with the fact that
$\widehat{F}_t^{\mathbb P,0}<+\infty$, $\P$-a.s. for every $\P\in
\mathcal{P}^\kappa_H$, implies that $\hat{a}_t \in D^1_{F_t}$
and $\nu^{\mathbb P}_t \in D^2_{F_t}$ $dt\times\P$-a.s., for all $\P
\in\mathcal{P}^\kappa_H$.
\end{Remark}

\subsection{The spaces and norms}\label{secspace}

We now define as in \cite{stz}, the spaces and norms which will be
needed for the formulation of the 2BSDEJs.

For $p\geq1$, $L^{p,\kappa}_H$ denotes the space of all $\mathcal
F_T$-measurable scalar r.v. $\xi$ with
\[
\llVert\xi\rrVert_{L^{p,\kappa}_H}^p:=\sup_{\mathbb{P}
\in\mathcal
{P}_H^\kappa}
\mathbb E^{\mathbb P} \bigl[\llvert\xi\rrvert^p \bigr]<+\infty.
\]

$\mathbb H^{p,\kappa}_H$ denotes the space of all $\mathbb
F^+$-predictable $\mathbb R^d$-valued processes $Z$ with
\[
\llVert Z\rrVert_{\mathbb H^{p,\kappa}_H}^p:=\sup_{\mathbb
{P} \in\mathcal
{P}_H^\kappa}
\mathbb E^{\mathbb P} \biggl[ \biggl(\int_0^T
\bigl\llvert\hat a_t^{1/2}Z_t \bigr\rrvert
^2\,dt \biggr)^{{p}/2} \biggr]<+\infty.
\]

$\mathbb D^{p,\kappa}_H$ denotes the space of all $\mathbb
F^+$-progressively measurable $\mathbb R$-valued processes $Y$ with
\[
\mathcal P^\kappa_H\mbox{-q.s. c\`adl\`ag paths, and }
\llVert Y\rrVert_{\mathbb D^{p,\kappa}_H}^p:=\sup_{\mathbb{P} \in
\mathcal
{P}_H^\kappa}
\mathbb E^{\mathbb P} \Bigl[\sup_{0\leq t\leq
T}\llvert
Y_t\rrvert^p \Bigr]<+\infty.
\]

$\mathbb J^{p,\kappa}_H$ denotes the space of all $\mathbb
F^+$-predictable functions $U$ with
\[
\llVert U\rrVert_{\mathbb J^{p,\kappa}_H}^p:=\sup_{\mathbb
{P} \in\mathcal
{P}_H^\kappa}
\mathbb E^{\mathbb P} \biggl[ \biggl(\int_0^T\!\!\int_{E} \bigl\llvert U_t(x) \bigr\rrvert
^2 \nu^{\mathbb P}_t(dx)\,dt \biggr)^{{p}/2}
\biggr]<+\infty.
\]

For each $\xi\in L^{1,\kappa}_H$, $\mathbb P\in\mathcal P^\kappa_H$
and $t \in[0,T]$, denote
\[
\mathbb E_t^{H,\mathbb P}[\xi]:=\mathop{\operatorname{ess}\operatorname
{sup}^{\mathbb P}}_{\mathbb P'\in
\mathcal P^\kappa_H(t^{+},\mathbb P)}\mathbb E^{\mathbb P'}_t[\xi],
\]
where
\[
\mathcal P^\kappa_H
\bigl(t^{+},\mathbb P \bigr):= \bigl\{ \mathbb P'\in
\mathcal P^\kappa_H\dvtx\mathbb P'=\mathbb P
\mbox{ on }\mathcal F_t^+ \bigr\}.
\]

Then we define for each $p\geq\kappa$,
\[
\mathbb L_H^{p,\kappa}:= \bigl\{\xi\in L^{p,\kappa}_H\dvtx
\llVert\xi\rrVert_{\mathbb L_H^{p,\kappa}}<+\infty\bigr\},
\]
where
\[
\llVert
\xi\rrVert_{\mathbb L_H^{p,\kappa}}^p:= \sup_{\mathbb P\in\mathcal
P^\kappa_H}\mathbb
E^{\mathbb P} \Bigl[\mathop{\operatorname{ess}\operatorname{sup}^{\mathbb P}}_{0\leq t\leq
T}\bigl(\mathbb E_t^{H,\mathbb
P} \bigl[\llvert\xi\rrvert
^\kappa\bigr] \bigr)^{{p}/{\kappa
}} \Bigr].
\]

%
\begin{Remark}
Except for $\mathbb L_H^{p,\kappa}$, the definitions of the previous
spaces are classic, but the second-order framework induces the presence
of an essential supremum over our family of probability measures. As
for $\mathbb L_H^{p,\kappa}$, it appears naturally in the a priori
estimates; we refer to \cite{stz} for more details.
\end{Remark}

Finally, we denote by $\UC_b(\Omega)$ the collection of all
bounded and uniformly continuous maps $\xi\dvtx\Omega\rightarrow\mathbb
R$ with respect to the Skorohod distance $d_S$, and we let
\begin{eqnarray}
\mathcal L^{p,\kappa}_H:=\mbox{the closure of }
\UC_b(\Omega)\mbox{ under the norm }\llVert\cdot\rrVert
_{\mathbb
L^{p,\kappa}_H}\nonumber
\\
\eqntext{\mbox{for every }1<\kappa\leq p.}
\end{eqnarray}

%
\begin{Remark}
In our accompanying paper \cite{kpz}, we will prove existence for
2BSDEJs for terminal conditions belonging to the space $\mathcal
L^{2,\kappa}_H$. We therefore think that it is important to give a few
examples of terminal conditions belonging to it. First of all, with
applications of 2BSDEJs to fully nonlinear PIDEs, we at least would
like functions of the form $f(B_T)$ to be in $\mathcal L^{2,\kappa
}_H$. But it is a well-known result that the application $\omega
\mapsto B_t(\omega)$ is continuous for the Skorohod topology for
Lebesgue almost every $t\in[0,T]$, including $t=0$ and $t=T$. Hence,
it is easy to see that for a Lipschitz function $f$, $f(B_t)\in
\mathcal L^{2,\kappa}_H$ for a.e. $t\in[0,T]$, including $t=0$ and
$t=T$. We also refer the reader to our accompanying paper \cite{kpz}
for more explanations and intuitions about this problem. Finally we
would like to mention that the recent results of \cite{nn2,nn3}, which
appeared during the revision of this paper, could be used to obtain
existence of a solution when $F=0$, but with a terminal condition $\xi
\in\L^{2,\kappa}_H$. It is a very interesting and difficult problem
to see whether their approach could be extended to general generators $F$.
\end{Remark}

For a given probability measure $\mathbb P\in\mathcal P^\kappa_H$,
the spaces $L^p(\mathbb P)$, $\mathbb D^p(\mathbb P)$, $\mathbb
H^p(\mathbb P)$ and $\mathbb J^p(\mathbb P)$ correspond to the above
spaces when the set of probability measures is only the singleton
$ \{\mathbb P \}$. Finally, we have
$\mathbb H^{p}_{\mathrm{loc}}(\mathbb P)$ denotes the space of all $\mathbb
F^+$-predictable $\mathbb R^d$-valued processes $Z$ with
\[
\biggl(\int_0^T \bigl\llvert\hat
a_t^{1/2}Z_t \bigr\rrvert^2\,dt
\biggr)^{{p}/2}<+\infty,\qquad \mathbb P\mbox{-a.s.}
\]

$\mathbb J^{p}_{\mathrm{loc}}(\mathbb P)$ denotes the space of all $\mathbb
F^+$-predictable functions $U$ with
\[
\biggl(\int_0^T\!\!\int_{E}
\bigl\llvert U_t(x) \bigr\rrvert^2 \nu
^{\mathbb P}_t(dx)\,dt \biggr)^{{p}/2}<+\infty,\qquad \mathbb P\mbox{-a.s.}
\]

\subsection{Formulation}

We shall consider the following $2$BSDEJ, for $0\leq t\leq T$ and
$\mathcal{P}^\kappa_H\mbox{-q.s.}$:
%
\begin{eqnarray}\label{2bsdej}
Y_t&=&\xi+\int_t^T
\widehat{F}^{\mathbb P}_s(Y_s,Z_s,U_s)\,ds
-\int_t^T Z_s\,dB^{\mathbb P,c}_s
\nonumber\\[-8pt]\\[-8pt]\nonumber
&&{} -
\int_t^T\!\!\int_{E}
U_s(x) \tilde{\mu}^{\mathbb
P}_B(dx,ds) +
K^{\mathbb P}_T-K^{\mathbb P}_t.
\end{eqnarray}

%
\begin{Definition}
We say $(Y,Z,U)\in\mathbb D^{2,\kappa}_H \times\mathbb H^{2,\kappa
}_H \times\mathbb J^{2,\kappa}_H$ is a solution to $2$BSDEJ (\ref
{2bsdej}) if:
\begin{itemize}
\item $Y_T=\xi$, $\mathcal{P}^\kappa_H$-q.s.
\item For all $\mathbb P \in\mathcal{P}^\kappa_H$ and
$0\leq t\leq T$, the process $K^{\mathbb P}$ defined below is
predictable and has nondecreasing paths $\mathbb P$-a.s.
%
\begin{eqnarray}\label{2bsdejK}
K_t^{\mathbb P}&:=&Y_0-Y_t - \int
_0^t\widehat{F}^{\mathbb
P}_s(Y_s,Z_s,U_s)\,ds+
\int_0^tZ_s\,dB^{\mathbb P,c}_s
\nonumber\\[-8pt]\\[-8pt]\nonumber
&&{} + \int_0^t\!\int_{E}
U_s(x) \tilde{\mu}^{\mathbb P}_B(dx,ds).
\end{eqnarray}

\item The family $ \{K^{\mathbb P}, \mathbb P \in
\mathcal P_H^\kappa\}$ satisfies the minimum condition
%
\begin{equation}
K_t^{\mathbb P}= \mathop{\operatorname{ess}\operatorname{inf}^{\mathbb P}}_{
\mathbb{P}' \in\mathcal
{P}_H(t^+,\mathbb{P})} \mathbb{E}_t^{\mathbb P'} \bigl[K_T^{\mathbb
{P}'}
\bigr],\qquad0\leq t\leq T,\mathbb P\mbox{-a.s.}, \forall
\mathbb P \in\mathcal P_H^\kappa. \label{2bsdejminK}
\end{equation}
\end{itemize}
Moreover if the family $ \{K^{\mathbb P}, \mathbb P \in\mathcal
P_H^\kappa\}$ can be aggregated into a universal process~$K$,
we call $(Y,Z,U,K)$ a solution of the $2$BSDEJ (\ref{2bsdej}).
\end{Definition}


Following \cite{stz}, in addition to Assumption \ref{assumphref}, we
will always assume the following:

%
\begin{Assumption}\label{assumph2ref}
\textup{(i)}~$\mathcal P_H^\kappa$ is not empty.

{(ii)} The process $F$ satisfies the following integrability
condition: 
%
\begin{equation}
\phi^{2,\kappa}_H:=\sup_{\mathbb P\in\mathcal P^\kappa_H}\mathbb
E^{\mathbb P} \biggl[\mathop{\operatorname{ess}\operatorname{sup}^{\mathbb P}}_{0\leq
t\leq T}
\biggl(\mathbb E_t^{H,\mathbb P} \biggl[\int^T_0
\bigl\llvert\Fh^{\mathbb
P,0}_s \bigr\rrvert^\kappa \,ds
\biggr] \biggr)^{{2}/{\kappa}} \biggr]<+\infty.
\end{equation}
\end{Assumption}

\subsection{Connection with standard BSDEJs}\label{connectionbsdes}

Let us assume that $H$ is linear in $\gamma$ and $\widetilde{v}$, in the
following sense:
%
\begin{equation}
H_t(y,z,u,\gamma,\widetilde{v}):=\frac{1}2\operatorname{Tr}
[I_d\gamma] + \int_E \widetilde{v}(e)
\nu^*(de) -f_t(y,z,u),\label{linearH}
\end{equation}
where $\nu^* \in\mathcal{N}$. We then have the following result:

%
\begin{Lemma}
If $H$ is of the form (\ref{linearH}), then $D^1_{F_t} = \{I_d\}$,
$D^2_{F_t}=\{\nu^*\}$ and
\[
F_t(\omega,y,z,u,a,\nu) = F_t \bigl(\omega,y,z,u,Id,
\nu^* \bigr) = f_t(y,z,u).
\]
\end{Lemma}

\begin{pf}
First notice that
\begin{eqnarray*}
&& H_t(\omega,y,z,u,\gamma,\widetilde{v})
\\
&&\qquad = \sup_{(a,\nu) \in
\mathbb{S}_d^{>0} \times\mathcal{N}}
\biggl\{ \frac{1}{2} \trace(a\gamma) + \int_0^T\!\!\int_E \widetilde{v}(e)\nu_s(\omega)
(ds,de) -\delta_{Id}(a) - \delta_{\nu^*}(\nu) \biggr\}
\\
&&\quad\qquad{}-f_t(y,z,u),
\end{eqnarray*}
where $\delta_A$ denotes the characteristic function of a subset $A$
in the convex analysis sense.

By definition of $F$, we get
\[
F_t(\omega,y,z,u,a,\nu)=f_t(y,z,u) + H^{**}(a,
\nu),
\]
where $H^{**}$ is the double Fenchel--Legendre transform of the function
\[
(a,\nu) \mapsto\delta_{Id}(a) + \delta_{\nu^*}(\nu),
\]
which is convex and lower-semicontinuous.

This then implies that
\[
F_t(\omega,y,z,u,a,\nu) =f_t(y,z,u) +
\delta_{Id}(a) + \delta_{\nu
^*}(\nu),
\]
which is the desired result.
\end{pf}

If we further assume that $\mathbb E^{\P_{\nu^*}} [\int_0^T\llvert
f_t(0,0,0)\rrvert^2\,dt ]<+\infty$, then
$\mathcal P^\kappa_H= \{
\P_{\nu^*} \}$ and the minimum condition on $K=K^{\P_{\nu
^*}}$ implies that $0=\mathbb E^{\P_{\nu^*}} [K_T ]$,
which means that $K\equiv0$, $\P_{\nu^*}$-a.s., and the 2BSDEJ is
reduced to a classical BSDEJ.

\subsection{Connection with $G$-expectations and $G$-L\'evy processes}\label{connectionglevy}

\subsubsection{Reminder on $G$-L\'evy processes}

In their recent paper, Hu and Peng \cite{hupengglevy} introduced a
new class of processes with independent and stationary increments,
called $G$-L\'evy processes. These processes are defined intrinsically,
that is, without making reference to any probability measure.\vadjust{\goodbreak}

Let $\widetilde{\Omega}$ be a given set, and let $\Hc$ be a linear
space of real valued functions defined on $\widetilde{\Omega}$,
containing the constants and such that $\llvert X\rrvert\in
\Hc$ if $X \in
\Hc$.
A sublinear expectation is a functional $\widehat{\E}\dvtx \Hc
\rightarrow
\R$ which is monotone increasing, constant preserving, sub-additive
and positively homogeneous. We\vspace*{1pt} refer to Definition 1.1 of \cite{peng2}
for more details. The triple $(\widetilde{\Omega}, \Hc, \widehat{\E})$
is called a sublinear expectation space.

%
\begin{Definition}
A $d$-dimensional c\`adl\`ag process $\{X_t, t \geq0\}$ defined on
a sublinear expectation space $(\widetilde{\Omega}, \Hc, \widehat
{\E
})$ is called a $G$-L\'evy process if:
\begin{longlist}[(iii)]
\item[(i)] $X_0 =0$.
\item[(ii)] $X$ has independent increments: $\forall s,t >0$,
the random variable $(X_{t+s} - X_t)$ is independent from
$(X_{t_1},\ldots, X_{t_n})$, for each $n\in\N$ and $0\leq t_1 < \cdots<t_n
\leq t$.
The notion of independence used here corresponds to Definition 3.10 in
\cite{peng2}.
\item[(iii)] $X$ has stationary increments: $\forall s,t >0$,
the distribution of $(X_{t+s} - X_t)$ does not depend on $t$. The
notion of distribution used here corresponds to the definition given in
Section~3 of \cite{peng2}.
\item[(iv)] For each $t \geq0$, there exists a decomposition
$X_t = X_t^c + X_t^d$, where $\{X_t^c, t \geq0\}$ is a continuous
process and $\{X_t^d, t \geq0\}$ is a pure jump process.
\item[(v)] $(X_t^c, X_t^d)$ is a $2d$-dimensional process
satisfying conditions (i), (ii) and (iii) of this definition and
\[
\lim_{t \to0^+} \frac{1}{t} \widehat{\E} \bigl( \bigl\llvert
X_t^c \bigr\rrvert^3 \bigr)=0,\qquad \widehat{\E}
\bigl( \bigl\llvert X_t^d \bigr\rrvert\bigr) \leq C t,  t
\geq0
\]
for a real constant $C$.
\end{longlist}
\end{Definition}

In \cite{hupengglevy}, Hu and Peng proved the following L\'
evy--Khintchine representation for $G$-L\'evy processes:

\begin{Theorem}[(Hu and Peng \cite{hupengglevy})] \label{glevykhintchine}
Let $\{X_t, t \geq0\}$ be a $G$-L\'evy process. Then for each
Lipschitz and bounded function $\varphi$, the function $u$ defined by
$u(t,x):= \widehat{\E} (\varphi[x+X_t] )$ is the unique
viscosity solution of the following partial integro-differential equation:
\begin{eqnarray*}
&& \partial_t u(t,x) -\sup_{(b,\alpha,\nu)\in\Uc} \biggl\{ \int
_E \bigl[u(t,x+z)-u(t,x) \bigr]\nu(dz)
\\
&&\hspace*{92pt}{} + \bigl\langle
Du(t,x),b \bigr\rangle+ \frac{1}{2}\operatorname{Tr} \bigl[D^2u(t,x)
\alpha\alpha^T \bigr] \biggr\}=0,
\end{eqnarray*}
where $\Uc$ is a subset of $\R^d \times\R^{d\times d} \times
\mathcal{M}_R^+$ satisfying
\[
\sup_{(b,\alpha,\nu)\in\Uc} \biggl\{ \int_{\R^d} \llvert z
\rrvert\nu(dz) + \llvert b\rrvert+ \operatorname{Tr} \bigl[\alpha
\alpha^T \bigr] \biggr\} < + \infty,
\]
and where $\mathcal{M}_R^+$ denotes the set of positive Radon measures
on $E$.
\end{Theorem}

Notice that Hu and Peng study the case of $G$-L\'evy processes with a
discontinuous part that is of finite variation.

\subsubsection{A connection with a particular 2BSDEJ}

In our framework, we know that $B^d$ is a purely discontinuous
semimartingale of finite variation under $\mathbb P_{0,F}$ if
\[
\int_0^T\!\!\int_{\llvert x\rrvert\leq1} \llvert
x\rrvert F_s(dx)\,ds < + \infty,\qquad \P_{0,F}\mbox{-a.s.}
\]
We give a function $H$ below, which is the natural candidate to
retrieve the example of $G$-L\'evy processes in our context. This link
will be made clear in our accompanying paper \cite{kpz}.\vspace*{1pt}

Let $\widetilde{\mathcal{N}}$ be any subset of $\mathcal{V}$ that is
convex and closed for the weak topology on $\mathcal{M}_R^+$. We define
\begin{eqnarray*}
&& H_t(\omega, \gamma, \widetilde{v})
\\
&&\qquad := \sup_{(a,\nu) \in\mathbb
{S}_d^{>0} \times\mathcal{N}}
\biggl\{ \frac{1}{2} \trace(a\gamma) + \int_0^T\!\!\int_E \widetilde{v}(e)\nu_s(de)\,ds -\delta
_{[a_1,a_2]}(a) - \delta_{\widetilde{\mathcal{N}}}(\nu) \biggr\}.
\end{eqnarray*}
Since $[a_1,a_2]$ and $\widetilde{\mathcal{N}}$ are closed convex spaces,
$F_t(\omega, a, \nu)$ is the double Fenchel--Legendre transform in
$(a,\nu)$ of the convex and lower semi-continuous function $(a,\nu)
\mapsto\delta_{[a_1,a_2]}(a) + \delta_{\widetilde{\mathcal
{N}}}(\nu)$
and then
\[
F_t(\omega, a, \nu) = \delta_{[a_1,a_2]}(a) +
\delta_{\widetilde
{\mathcal{N}}}(\nu).
\]
In \cite{kpz}, we prove that the 2BSDEJs are connected to a class of
fully nonlinear partial integro-differential equations. With this
particular function $H$ and its transform $F$, the PIDE we find is the
one given in Theorem \ref{glevykhintchine}.
If moreover $H_t(\omega, \gamma, \widetilde{v}) = H_t(\omega,
\tilde{v})$ is independent of $\gamma$, and $\widetilde{\mathcal
{N}} = \{
\lambda\delta_{\{1\}}, \lambda_1\leq\lambda\leq\lambda_2 \}$
(which is convex and closed and where $\delta_{\{1\}}$ is a Dirac mass
at the point $1$), then $F_t$ is independent of $a$, and we obtain a
2BSDEJ giving a representation of the $G$-Poisson process.

\section{Uniqueness result}\label{section2}

In this section, we address the question of uniqueness of a solution to
a 2BSDEJ. We follow the intuition provided in the \hyperref[sec1]{Introduction} and
write the solution to a 2BSDEJ as a supremum in some sense of solutions
to classical BSDEJs.

\subsection{Representation of the solution}
We have the following, which is similar to Theorem 4.4 of \cite{stz}:


\begin{Theorem}\label{uniqueref}
Let Assumptions \ref{assumphref} and \ref{assumph2ref} hold. Assume
$\xi\in\mathbb{L}^{2,\kappa}_H$ and that $(Y,Z,U)$ is a solution to
the $2$BSDEJ (\ref{2bsdej}). Then, for any $\mathbb{P}\in\mathcal
{P}^\kappa_H$ and $0\leq t_1< t_2\leq T$,
%
%
\begin{eqnarray}\label{representationref}
Y_{t_1}&=&\mathop{\operatorname{ess}\operatorname
{sup}^\mathbb{P}}_{\mathbb
{P}'\in\mathcal{P}^\kappa
_H(t_1^+,\mathbb{P})}y_{t_1}^{\mathbb{P}'}(t_2,Y_{t_2}),
\qquad\mathbb{P}\mbox{-a.s.},
\end{eqnarray}
where, for any $\mathbb{P}\in\mathcal{P}^\kappa_H$, $\mathbb
{F}^+$-stopping time $\tau$, and $\mathcal{F}^+_{\tau} $-measurable
random variable $\xi\in\mathbb{L}^2({\mathbb P})$, $(y^{\mathbb
{P}}(\tau,\xi),z^{\mathbb{P}}(\tau,\xi))$ denotes the solution to
the following standard BSDEJ on $0\leq t\leq\tau$
%
\begin{eqnarray}\label{bsdej}
y^{\mathbb{P}}_t &=& \xi+ \int_t^{\tau}
\widehat{F}^{\mathbb
{P}}_s \bigl(y^{\mathbb{P}}_s,z^{\mathbb{P}}_s,u^{\mathbb{P}}_s
\bigr)\,ds-\int_t^{\tau}z^{\mathbb{P}}_s\,dB^{\mathbb{P},c}_s
\nonumber\\[-8pt]\\[-8pt]\nonumber
&&{} - \int_t^{\tau}\!\int_{E}
u^{\mathbb{P}}_s(x) \tilde{\mu}^{\mathbb{P}}_B(dx,ds),
\qquad\mathbb P\mbox{-a.s.}
\end{eqnarray}

Consequently, the $2$BSDEJ (\ref{2bsdej}) has at most one solution in $
\mathbb D^{2,\kappa}_H\times\mathbb H^{2,\kappa}_H\times\mathbb
J^{2,\kappa}_H$.
\end{Theorem}

%
\begin{Remark}\label{rembsde}
We first emphasize that existence and uniqueness results for the
standard BSDEJs (\ref{bsdej}) are not given directly by the existing
literature, since the compensator of the counting measure associated to
the jumps of $B$ is not deterministic. However, since all the
probability measures we consider satisfy the martingale representation
property and the Blumenthal 0--1 law, it is clear that we can
straightforwardly generalize the proof of existence and uniqueness of
Tang and Li \cite{tangli}; see also \cite{bech} and \cite{crep} for
related results. Furthermore, the usual a priori estimates and
comparison theorems will also hold.
\end{Remark}

%

Before giving the proof of the above theorem, we first state the
following lemma which is a generalization of the usual comparison
theorem proved by Royer; see Theorem 2.5 in \cite{roy}. Its proof is
a straightforward generalization so we omit it.

%
\begin{Lemma}\label{lemmacompref} Let $\mathbb P\in\mathcal P^\kappa
_H$. We consider two generators $f^1$ and $f^2$ satisfying Assumption
$H_{comp}$ in \cite{roy} (which\vspace*{1pt} is a consequence of our more
restrictive assumptions). Given two nondecreasing processes $k^1$ and
$k^2$, let $\xi^1$ and $\xi^2$ be two terminal conditions for the
following BSDEJs for $i=1,2$,
\begin{eqnarray*}
y^i_{t} &=& \xi^i + \int
_{t}^{T} f^i_s
\bigl(y^i_s,z^i_s,u^i_s
\bigr)\,ds - \int_{t}^{T} z^i_s\,dB_s
\\
&&{} - \int_{t}^{T}\!\!\int_E
u^i_s(x)\widetilde\mu^\P(dx,ds) +
k^i_{T}-k^i_{t},\qquad \mathbb P\mbox{-a.s.}
\end{eqnarray*}

Denote by $(y^1, z^1,u^1)$ and $(y^2, z^2,u^2)$ the respective
solutions. If $\xi^1 \leq\xi^2$, $k^1-k^2$ is nonincreasing and
$f^1(t, y^1_t, z^1_t,u^1_t ) \leq f^2(t, y^1_t, z^1_t,u^1_t )$,
then $\forall t \in[0, T ], Y^1_t \leq Y^2_t$.
\end{Lemma}

\begin{pf*}{Proof of Theorem \ref{uniqueref}}
The proof follows the lines of the proof of Theorem 4.4 in \cite
{stz}. First of all, if representation (\ref{representationref}) holds,
then $Y$ is uniquely defined. Moreover, since we have that
\[
d[Y,B]^c_t=Z_td[B,B]^c_t=
\hat a_tZ_t\,dt,\qquad\mathcal{P}^\kappa_H\mbox{-q.s.},
\]
$Z$ is also uniquely defined.

Then, since for any $\P\in\Pc^\kappa_H$, $B$ only has $\P$-totally
inaccessible jump times, we know that $B$ and $K^\P$ never jump at the
same time, $\P$-a.s. We deduce that
%
\begin{equation}
\label{eqjump} \Delta[Y,B]_t=U_t(\Delta B_t)
\Delta B_t{\mathbf1}_{\Delta B_t\neq
0},\qquad\mathcal{P}^\kappa_H
\mbox{-q.s.}
\end{equation}

We can then define
\[
\widetilde U_t(x):=U_t(x){\mathbf1}_{\Delta B_t=x}.
\]

Then $\widetilde U$ is actually equal to $U$, $dt\times\nu_t^\P
(dx)$, for any $\P\in\mathcal P^\kappa_H$. Using this version
instead, and still denoting it $U$ for simplicity, we deduce that $U$
is uniquely defined by (\ref{eqjump}). Then the uniqueness of the
process $K^{\mathbb P}$ is immediate. Let us now proceed with the proof
of (\ref{representationref}):
\begin{longlist}[(iii)]
\item[(i)] Fix $0\leq t_1<t_2\leq T$ and $\mathbb P\in\mathcal P^\kappa
_H$. For any $\mathbb P'\in\mathcal P^\kappa_H(t_1^+,\mathbb P)$ and
$t_1\leq t\leq t_2$, we have
\begin{eqnarray*}
Y_{t} &=& Y_{t_2} + \int_{t}^{t_2}
\widehat{F}^{\mathbb
P'}_s(Y_s,Z_s,U_s)\,ds
- \int_{t}^{t_2} Z_s\,dB^{\mathbb P',c}_s
\\
&&{} - \int_{t}^{t_2}\!\!\int_E
U_s(x)\widetilde\mu^{\mathbb P'}_B(dx,ds)
+ K^{\mathbb P'}_{t_2}-K^{\mathbb P'}_{t},\qquad
\mathbb P'\mbox{-a.s.}
\end{eqnarray*}

With Assumption \ref{assumphref}, we can apply the above Lemma \ref
{lemmacompref} under $\mathbb P'$ to obtain that $Y_{t_1}\geq
y_{t_1}^{\mathbb P'}(t_2,Y_{t_2}), \mathbb P'\mbox{-a.s.}$
Since $\mathbb P'=\mathbb P$ on $\mathcal F_{t_1}^+$, we get
$Y_{t_1}\geq y_{t_1}^{\mathbb P'}(t_2,Y_{t_2})$, $\mathbb P$-a.s. and thus
\[
Y_{t_1}\geq \mathop{\operatorname{ess}\operatorname
{sup}^\mathbb{P}}_{\mathbb{P}'\in
\mathcal{P}^\kappa_H(t_1^+,\mathbb{P})}y_{t_1}^{\mathbb
{P}'}(t_2,Y_{t_2}),
\qquad\mathbb{P}\mbox{-a.s.}
\]

\item[(ii)] We now prove the reverse inequality. Fix $\mathbb P\in
\mathcal P^\kappa_H$. We will show in (iii) below that
\[
C_{t_1}^{\mathbb P}:=\mathop{\operatorname{ess}\operatorname{sup}^\mathbb{P}}_{\mathbb{P}'\in\mathcal{P}^\kappa
_H(t_1^+,\mathbb{P})}\mathbb E_{t_1}^{\mathbb{P}'}
\bigl[ \bigl(K^{\mathbb{P}'}_{t_2}-K^{\mathbb
{P}'}_{t_1}
\bigr)^2 \bigr]<+\infty,\qquad\mathbb P\mbox{-a.s.}
\]

For every $\mathbb P'\in\mathcal P^\kappa_H(t^+,\mathbb P)$, denote
\begin{eqnarray*}
\delta Y&:=&Y-y^{\mathbb P'}(t_2,Y_{t_2}),\qquad\delta
Z:=Z-z^{\mathbb
P'}(t_2,Y_{t_2})\quad\mbox{and}
\\
\delta U&:=&U-u^{\mathbb P'}(t_2,Y_{t_2}).
\end{eqnarray*}

By the Lipschitz Assumption \ref{assumphref}{(iii)},
there exist two bounded processes $\lambda$ and $\eta$ such that for
all $t_1\leq t\leq t_2$,
\begin{eqnarray*}
\delta Y_t&=& \int_t^{t_2} \bigl(
\lambda_s\delta Y_s+\eta_s\hat
{a}_s^{1/2}\delta Z_s \bigr)\,ds
\\
&&{} +\int
_t^{t_2} \bigl(\widehat{F}^{\mathbb{P}'}_s
\bigl(y^{\mathbb P'}_s,z^{\mathbb P'}_s,U_s
\bigr)-\widehat{F}^{\mathbb{P}'}_s \bigl(y^{\mathbb P'}_s,z^{\mathbb
P'}_s,u^{\mathbb
P'}_s
\bigr) \bigr)\,ds
\\
&&{}-\int_t^{t_2}\delta Z_s\,dB^{\mathbb{P}',c}_s-
\int_t^{t_2}\!\!\int_E \delta
U_s(x)\widetilde\mu^{\mathbb{P}'}_B(dx,ds)+
K^{\mathbb
{P}'}_{t_2}- K^{\mathbb{P}'}_{t},\qquad\mathbb
P'\mbox{-a.s.}
\end{eqnarray*}

Define for $t_1\leq t\leq t_2$ the following processes:
\[
N^{\mathbb{P}'}_t:=\int_{t_1}^t
\eta_s\hat{a}_s^{-1/2}\,dB^{\mathbb
{P}',c}_s+
\int_{t_1}^t\int_{E}\gamma
_s(x)\widetilde\mu^{\mathbb{P}'}_B(ds, dx)
\]
and
\[
M^{\mathbb{P}'}_t:=\exp\biggl(\int_{t_1}^t
\lambda_s\,ds \biggr)\mathcal E \bigl(N^{\mathbb{P}'}
\bigr)_t,
\]
where $\mathcal E(N^{\mathbb{P}'})_t$ denotes the Dol\'eans--Dade
exponential martingale of $N^{\mathbb{P}'}_t$.

By the boundedness of $\lambda$ and $\eta$ and the assumption on
$\gamma$ in Assumption \textup{\ref{assumphref}(iv)}, we know that
$M$ has moments (positive or negative) of any order; see \cite
{lepinmemin2} for the positive moments and Lemma \ref{expomart2} in
the \hyperref[app]{Appendix} for the negative ones. Thus we have for $p\geq1$,
%
\begin{equation}
\label{trucM} \mathbb E_{t_1}^{\mathbb P'} \Bigl[\sup
_{t_1\leq t\leq t_2} \bigl(M^{\mathbb
{P}'}_t
\bigr)^p+ \sup_{t_1\leq t\leq t_2} \bigl(M^{\mathbb{P}'}_t
\bigr)^{-p} \Bigr]\leq C_p, \qquad\mathbb
P' \mbox{-a.s.}
\end{equation}

Then, by It\^o's formula, we obtain
\begin{eqnarray*}
&& d \bigl(M^{\mathbb P'}_t\delta Y_{t} \bigr)
\\
&&\qquad =
M^{\mathbb P'}_{t-}\,d(\delta Y_{t})+\delta
Y_{t-} \,d M^{\mathbb P'}_t+d \bigl[M^{\mathbb P'},
\delta Y \bigr]_t
\\
&&\qquad =  M^{\mathbb P'}_{t-} \biggl[ \bigl(-\lambda_t \delta
Y_t-\eta_t\hat{a}_t^{1/2}\delta
Z_t-\widehat{F}^{\mathbb P'}_t \bigl(y^{\mathbb
P'}_t,z^{\mathbb P'}_t,U_t
\bigr)+\widehat{F}^{\mathbb P'}_t \bigl(y^{\mathbb
P'}_t,z^{\mathbb P'}_t,u^{\mathbb P'}_t
\bigr) \bigr)\,dt
\\
&&\hspace*{94pt}{}+\delta Z_t\,dB^{\mathbb P',c}_t+\int
_E \bigl(\delta U_t(x)+\gamma_t(x)
\delta U_t(x) \bigr)\widetilde\mu^{\mathbb
P',c}_B(dx,dt)
\biggr]
\\
&&\quad\qquad{} +\delta Y_{t-}M^{\mathbb P'}_{t-} \biggl(
\lambda_t\,dt+\eta_t\hat{a}_t^{-1/2}\,dB^{\mathbb P',c}_t+
\int_E \gamma_t(x)\widetilde\mu
^{\mathbb P'}_B(dx,dt) \biggr)
\\
&&\quad\qquad{} +M^{\mathbb P'}_{t} \biggl(\eta_t
\hat{a}_t^{1/2}\delta Z_t\,dt+ \int
_E \gamma_t(x)\delta U_t(x)
\nu^{\mathbb P'}_t(dx)\,dt \biggr)-M^{\mathbb P'}_{t-}\,dK^{\mathbb P'}_t.
\end{eqnarray*}

Thus, by Assumption \textup{\ref{assumphref}(iv)}, we have
\begin{eqnarray*}
\delta Y_{t_1}&\leq& -\int^{t_2}_{t_1}M^{\mathbb P'}_{s}
\bigl(\delta Z_s+ \delta Y_{s}\eta_s
\hat{a}_s^{-1/2} \bigr)\,dB^{\mathbb
P',c}_s+
\int^{t_2}_{t_1}M^{\mathbb P'}_{s-}\,dK^{\mathbb P'}_s
\\
&&{}-\int^{t_2}_{t_1}M^{\mathbb P'}_{s-}
\int_E \bigl(\delta U_s(x)+\delta
Y_{s}\gamma_s(x)+\gamma_s(x)\delta
U_s(x) \bigr)\widetilde\mu^{\mathbb P'}_B(dx,ds).
\end{eqnarray*}

By taking conditional expectation, we obtain
%
\begin{equation}
\label{brubru} \delta Y_{t_1}\leq\mathbb E_{t_1}^{\mathbb P'}
\biggl[\int_{t_1}^{t_2}M^{\mathbb P'}_{t-}
\,dK^{\mathbb P'}_t \biggr].
\end{equation}

Applying the H\"older inequality, we can now write
\begin{eqnarray*}
\delta Y_{t_1}&\leq&\mathbb E_{t_1}^{\mathbb P'} \Bigl[\sup
_{t_1\leq
t\leq t_2} \bigl(M^{\mathbb P'}_t \bigr) \bigl(
K_{t_2}^{\mathbb P'}- K_{t_1}^{\mathbb P'} \bigr) \Bigr]
\\
&\leq&\Bigl(\mathbb E_{t_1}^{\mathbb P'} \Bigl[\sup
_{t_1\leq t\leq
t_2} \bigl(M^{\mathbb P'}_t
\bigr)^3 \Bigr] \Bigr)^{1/3} \bigl(\mathbb
E_{t_1}^{\mathbb P'} \bigl[ \bigl( K_{t_2}^{\mathbb P'}-
K_{t_1}^{\mathbb P'} \bigr)^{3/2} \bigr]
\bigr)^{2/3}
\\
&\leq& C \bigl(C_{t_1}^\mathbb P \bigr)^{1/3} \bigl(
\mathbb E_{t_1}^{\mathbb
P'} \bigl[ K^{\mathbb P'}_{t_2}-
K^{\mathbb P'}_{t_1} \bigr] \bigr)^{1/3},\qquad\mathbb P
\mbox{-a.s.}
\end{eqnarray*}

Taking the essential infimum on both sides completes the proof.\vspace*{1pt}

\item[(iii)] It remains to show that the estimate for $C_{t_1}^{\mathbb
P}$ holds. But by definition, and the Lipschitz assumption on $F$, we
clearly have
%
\begin{eqnarray}\label{estimK}
&& \sup_{\mathbb{P}\in\mathcal{P}^\kappa_H} \mathbb
E^{\mathbb{P}} \bigl[
\bigl(K^{\mathbb P}_{t_2}-K^{\mathbb P}_{t_1}
\bigr)^2 \bigr]
\nonumber\\[-8pt]\\[-8pt]\nonumber
&&\qquad \leq C \bigl(\llVert Y\rrVert^2_{\mathbb
D^{2,\kappa}_H}+
\llVert Z\rrVert^{2}_{\H^{2,\kappa
}_H}+\llVert U\rrVert
^{2}_{\mathbb J^{2,\kappa}_H}+\phi^{2,\kappa}_H \bigr)<+
\infty,
\end{eqnarray}
since the last term on the right-hand side is finite thanks to the
integrability assumed on $\xi$ and $ F$. We then use the definition of
the essential supremum (see Neveu \cite{neveu}, e.g.) to have the
following equality:
%
\begin{equation}
\label{trucK}
\qquad\mathop{\operatorname{ess}\operatorname
{sup}^\mathbb
{P}}_{\mathbb{P}'\in\mathcal
{P}^\kappa_H(t_1^+,\mathbb{P})} \mathbb E^{\mathbb{P}'}_{t_1} \bigl[
\bigl(K^{\mathbb P'}_{t_2}-K^{\mathbb P'}_{t_1}
\bigr)^2 \bigr] =\sup_{n\geq1} \mathbb
E^{\mathbb{P}_{n}}_{t_1} \bigl[ \bigl(K^{\mathbb P_n}_{t_2}-K^{\mathbb
P_n}_{t_1}
\bigr)^2 \bigr],\qquad\mathbb{P}\mbox{-a.s.}
\end{equation}
for some sequence $ (\mathbb{P}_{n})_{n\geq1}\subset\mathcal
{P}^\kappa_H(t_1^+,\mathbb{P})$. Moreover, in Lemma \ref{upward} of
the \hyperref[app]{Appendix}, it is proved that the set $\mathcal{P}^\kappa
_H(t_1^+,\mathbb{P})$ is upward directed which means that for any
$\mathbb{P}'_1, \mathbb{P}'_2\in\mathcal{P}^\kappa
_H(t_1^+,\mathbb{P})$, there exists $\mathbb{P}'\in\mathcal
{P}^\kappa_H(t_1^+,\mathbb{P})$ such that
\begin{eqnarray*}
\mathbb E^{\mathbb{P}'}_{t_1} \bigl[ \bigl(K^{\mathbb
{P}'}_{t_2}-K^{\mathbb{P}'}_{t_1}
\bigr)^2 \bigr]&=&\max\bigl\{ \mathbb E^{\mathbb{P}'_1}_{t_1}
\bigl[ \bigl(K^{\mathbb
{P}'_1}_{t_2}-K^{\mathbb{P}'_1}_{t_1}
\bigr)^2 \bigr], \mathbb E^{\mathbb{P}'_2}_{t_1} \bigl[
\bigl(K^{\mathbb
{P}'_2}_{t_2}-K^{\mathbb{P}'_2}_{t_1}
\bigr)^2 \bigr] \bigr\}.
\end{eqnarray*}

Hence, by using a subsequence if necessary, we can rewrite (\ref
{trucK}) as
\begin{eqnarray*}
\mathop{\operatorname{ess}\operatorname{sup}^\mathbb
{P}}_{\mathbb{P}'\in\mathcal
{P}^\kappa_H(t_1^+,\mathbb{P})} \mathbb E^{\mathbb{P}'}_{t_1} \bigl[
\bigl(K^{\mathbb{P}'}_{t_2}-K^{\mathbb{P}'}_{t_1}
\bigr)^2 \bigr]&=& \mathop{\lim}_{n\rightarrow\infty}
\uparrow\mathbb E^{\mathbb{P}_{n}}_{t_1} \bigl[ \bigl(K^{\mathbb
{P}_{n}}_{t_2}-K^{\mathbb{P}_{n}}_{t_1}
\bigr)^2 \bigr],\qquad\mathbb{P}\mbox{-a.s.}
\end{eqnarray*}

With (\ref{estimK}), we can then complete the proof exactly as in the
proof of Theorem~4.4 in \cite{stz}.\quad\qed
\end{longlist}\noqed
\end{pf*}

Finally, the comparison theorem below follows easily from the classical
one for BSDEJs (see, e.g., Theorem 2.5 in \cite{roy}) and the
representation (\ref{representationref}).

%
\begin{Theorem}
Let $(Y,Z,U)$ and $(Y',Z',U')$ be the solutions of $2$BSDEJs with
terminal conditions $\xi$ and $\xi'$, generators $\widehat F$ and
$\widehat F'$, respectively (with\vspace*{1pt} the corresponding function $H$ and
$H'$), and let $(y^\mathbb P,z^\mathbb P,u^\mathbb P)$ and
$(y'^{\mathbb P},z'^{\mathbb P},u'^{\mathbb P})$ the solutions of the
associated BSDEJs. Assume that they both verify our Assumptions \ref
{assumphref} and \ref{assumph2ref} and that we have:
\begin{itemize}
\item$\mathcal P^\kappa_H\subset\mathcal P^\kappa_{H'}$;\vspace*{1pt}
\item$\xi\leq\xi'$, $\mathcal P^\kappa_H$-q.s.;\vspace*{1pt}
\item$\widehat F^\P_t(y'^{\mathbb P}_t,z'^{\mathbb P}_t,u'^{\mathbb
P}_t)\leq\widehat F^{\prime\P}_t(y'^{\mathbb P}_t,z'^{\mathbb
P}_t,u'^{\mathbb P}_t)$,
$\mathbb P$-a.s., for all $\mathbb P\in\mathcal P^\kappa_H$.
\end{itemize}
Then $Y\leq Y'$, $\mathcal P^\kappa_H$-q.s.
\end{Theorem}

\subsection{A priori estimates}

We conclude this section by showing some a priori estimates which will
be useful to obtain the existence of a solution in \cite{kpz}.

%
\begin{Theorem}\label{estimatesref}
Let Assumptions \ref{assumphref} and \ref{assumph2ref} hold. Assume
$\xi\in\mathbb L^{2,\kappa}_H$ and $(Y,Z,U)\in\mathbb D^{2,\kappa
}_H\times\mathbb H^{2,\kappa}_H\times\mathbb J^{2,\kappa}_H$ is a
solution to the 2BSDEJ (\ref{2bsdej}). Let $ \{(y^\mathbb
P,z^\mathbb P,u^\mathbb P) \}_{\mathbb P\in\mathcal P^\kappa
_H}$ be the\vspace*{1pt} solutions of the corresponding BSDEJs (\ref{bsdej}). Then
there exists a constant $C_\kappa$ such that
\begin{eqnarray*}
\llVert Y\rrVert^2_{\mathbb D^{2,\kappa}_H}+\llVert Z\rrVert
^2_{\mathbb H^{2,\kappa
}_H}+\llVert U\rrVert^2_{\mathbb J^{2,\kappa}_H}+
\sup_{\mathbb P\in\mathcal
P^\kappa_H}\mathbb E^{\mathbb{P}} \bigl[ \bigl\llvert
K_T^\P\bigr\rrvert^2 \bigr] &\leq&
C_{\kappa} \bigl(\llVert\xi\rrVert^2_{\mathbb
L^{2,\kappa}_H}+\phi
^{2,\kappa}_H \bigr),
\\
\sup_{\mathbb P\in\mathcal P^\kappa_H} \bigl\{ \bigl\llVert y^\mathbb P
\bigr
\rrVert^2_{\mathbb D^{2}(\mathbb P)}+ \bigl\llVert z^\mathbb P \bigr
\rrVert^2_{\mathbb
H^{2}(\mathbb P)}+ \bigl\llVert u^\mathbb P \bigr
\rrVert^2_{\mathbb
J^{2}(\mathbb P)} \bigr\}&\leq& C_{\kappa} \bigl(
\llVert\xi\rrVert^2_{\mathbb
L^{2,\kappa}_H}+\phi^{2,\kappa}_H
\bigr).
\end{eqnarray*}
\end{Theorem}

\begin{pf}
As in the proof of the representation formula in Theorem \ref
{representationref}, the Lipschitz Assumption \ref{assumphref}(iii)
implies that there exist two bounded processes $\lambda$ and
$\eta$ such that for all $t$, and $\mathbb P$-a.s.,
\begin{eqnarray*}
y_t^\mathbb P&=& \xi+ \int_t^{T}
\bigl(\lambda_sy_s^\mathbb P +\eta
_s\hat{a}_s^{1/2}z_s^\mathbb
P+\widehat{F}^\mathbb P_s \bigl(0,0,u_s^\mathbb
P \bigr) \bigr)\,ds
\\
&&{} -\int_t^{T}z_s^\mathbb
P\,dB^{\mathbb
P,c}_s-\int_t^{T}\!\!\int_E u_s^\mathbb P(x)\widetilde
\mu^{\mathbb P}_B(dx,ds).
\end{eqnarray*}

Define the following processes:
\[
N^{\mathbb P}_t:=\int_{t}^T
\eta_s\hat{a}_s^{-1/2}\,dB^{\mathbb
P,c}_s+ \int_{t}^T\!\!\int_{E}
\gamma_s(x)\widetilde\mu^{\mathbb
P}_B(dx,ds)
\]
and
\[
M_t:=\exp\biggl(\int_{t}^T
\lambda_s\,ds \biggr)\mathcal E \bigl(N^{\mathbb P}
\bigr)_t,
\]
where $\mathcal E(N^{\mathbb P})_t$ denotes the Dol\'eans--Dade
exponential martingale of $N^{\mathbb P}_t$. Then by applying It\^o's
formula to $M^{\mathbb P}_t y_t^\mathbb P$, we obtain
\[
y_t^\mathbb P = \mathbb E_t^\mathbb P
\biggl[M^{\mathbb P}_T\xi+ \int^T_t
M^{\mathbb P}_{s}\widehat{F}^{\mathbb P}_s
\bigl(0,0,u_s^\mathbb P \bigr)\,ds-\int^T_t\!\!
\int_E M^{\mathbb P}_{s}
\gamma_s(x)u_s^\mathbb P (x)
\nu^{\mathbb P}_s(dx)\,ds \biggr].
\]

Finally with Assumption \textup{\ref{assumphref}(iv)}, the H\"older
inequality and the inequality (\ref{trucM}), we conclude that there
exists a constant $C_\kappa$ depending only on $\kappa$, $T$ and the
Lipschitz constant of $F$, such that for all $\mathbb P$,
%
\begin{equation}
\label{estimref} \bigl\llvert y_t^\mathbb P \bigr\rrvert\leq
C_\kappa\mathbb E_t^\mathbb P \biggl[\llvert\xi
\rrvert^\kappa+\int_t^T \bigl\llvert
\widehat F^{\mathbb P,0}_s \bigr\rrvert^\kappa \,ds
\biggr]^{1/\kappa}.
\end{equation}

This immediately provides the estimate for $y^\mathbb P$. Now by
definition of our norms, we get from (\ref{estimref}) and
representation formula (\ref{representationref}) that
%
\begin{equation}
\label{eqy} \llVert Y\rrVert_{\mathbb D^{2,\kappa}_H}^2\leq
C_\kappa\bigl(\llVert\xi\rrVert^2_{\mathbb L^{2,\kappa}_H}+\phi
^{2,\kappa}_H \bigr).
\end{equation}

Now apply It\^o's formula to $\llvert Y\rrvert^2$ under each
$\mathbb P\in
\mathcal P^\kappa_H$. We get as usual for every $\epsilon>0$
\begin{eqnarray*}
&& \llvert Y_0\rrvert^2  +\int_0^T
\bigl\llvert\hat a_t^{1/2}Z_t \bigr\rrvert
^2\,dt+\int_0^T\!\!\int
_E \bigl\llvert U_t(x) \bigr\rrvert
^2\nu^{\mathbb P}_t(dx)\,dt
\\
&&\qquad = \llvert\xi\rrvert^2+2\int_0^TY_t
\widehat F^{\mathbb
P}_t(Y_t,Z_t,U_t)\,dt+2
\int_0^T Y_{t^-}\,dK^{\mathbb P}_t
\\
&&\quad\qquad{} -2\int_0^T Y_{t}Z_t\,dB^{\mathbb P,c}_t-
\int_0^T\!\!\int_E \bigl(
\bigl\llvert U_t(x) \bigr\rrvert^2+2Y_{t^-}U_t(x)
\bigr)\widetilde\mu^{\mathbb P}_B(dx,dt)
\\
&&\qquad \leq 2\int_0^T\llvert Y_t\rrvert
\bigl\llvert\widehat F^{\mathbb
P}_t(Y_t,Z_t,U_t)
\bigr\rrvert\, dt+2 \sup_{0\leq t\leq T} \llvert Y_{t}\rrvert
K^{\mathbb
P}_T
\\
&&\quad\qquad{}  -2\int_0^T Y_{t}Z_t\,dB^{\mathbb P,c}_t-
\int_0^T\!\!\int_E \bigl(
\bigl\llvert U_t(x) \bigr\rrvert^2+2Y_{t^-}U_t(x)
\bigr)\widetilde\mu^{\mathbb P}_B(dx,dt).
\end{eqnarray*}
By our assumptions on $F$, we have
\begin{eqnarray*}
&& \bigl\llvert\widehat F^{\mathbb P}_t(Y_t,Z_t,U_t)
\bigr\rrvert
\\
&&\qquad \leq C \biggl(\llvert Y_t\rrvert+ \bigl\llvert\hat
a_t^{1/2}Z_t \bigr\rrvert+ \bigl\llvert
\widehat F^{\mathbb P,0}_t \bigr\rrvert+ \biggl(\int
_E \bigl\llvert U_t(x) \bigr\rrvert
^2\nu^{\mathbb P}_t(dx) \biggr)^{1/2}
\biggr).
\end{eqnarray*}

With the usual inequality $2ab\leq\frac{1}{\epsilon} a^2+\epsilon
b^2, \forall\epsilon>0 $, we obtain
%
\begin{eqnarray}\label{grahou}
\nonumber
&&\mathbb E^\mathbb P \biggl[\int_0^T
\bigl\llvert\hat a_t^{1/2}Z_t \bigr\rrvert
^2\,dt+\int_0^T\!\!\int
_E \bigl\llvert U_t(x) \bigr\rrvert
^2\nu^{\mathbb
P}_t(dx)\,dt \biggr]
\\
\nonumber
&&\qquad \leq C\mathbb E^\mathbb P \biggl[\llvert\xi\rrvert
^2+\int_0^T\llvert Y_t
\rrvert\biggl( \bigl\llvert\widehat F^{\mathbb P,0}_t \bigr\rrvert
+ \llvert Y_t\rrvert+ \bigl\llvert\hat a_t^{1/2}Z_t \bigr\rrvert
\\
&&\hspace*{129pt}{} + \biggl(\int_E \bigl\llvert
U_t(x) \bigr\rrvert^2\nu^{\mathbb P}_t(dx)
\biggr)^{1/2} \biggr)\,dt \biggr]
\nonumber\\[-8pt]\\[-8pt]\nonumber
&&\quad\qquad{} +\mathbb E^\mathbb P \biggl[\int_0^T
\llvert Y_{t^-}\rrvert \,dK^{\mathbb P}_t \biggr]
\\
\nonumber &&\qquad \leq C \biggl(\llVert\xi\rrVert_{\mathbb
L^{2,\kappa}_H}+\mathbb
E^\mathbb P \biggl[ \biggl(1+\frac{C}\varepsilon\biggr)\sup
_{0\leq
t\leq
T}\llvert Y_t\rrvert^2+ \biggl(
\int_0^T \bigl\llvert\widehat
F^{\mathbb P,0}_t \bigr\rrvert \,dt \biggr)^2 \biggr]
\biggr)
\\
&&\quad\qquad{} +\epsilon\mathbb E^\mathbb P \biggl[\int_0^T
\bigl\llvert\hat a_t^{1/2}Z_t \bigr\rrvert
^2\,dt+\int_0^T\!\!\int
_E \bigl\llvert U_t(x) \bigr\rrvert
^2\nu^{\mathbb P}_t(dx)\,dt+ \bigl\llvert
K^{\mathbb P}_T \bigr\rrvert^2 \biggr].\nonumber
\end{eqnarray}

Then by definition of our $2$BSDEJ, we easily have
%
\begin{eqnarray}\label{eqkk}
\nonumber
\qquad \mathbb E^\mathbb P \bigl[ \bigl\llvert K^{\mathbb P}_T
\bigr\rrvert^2 \bigr]&\leq& C_0\mathbb E^\mathbb
P \biggl[\llvert\xi\rrvert^2+\sup_{0\leq t\leq
T} \llvert
Y_t\rrvert^2+\int_0^T
\bigl\llvert\hat a_t^{1/2}Z_t \bigr\rrvert
^2\,dt
\nonumber\\[-8pt]\\[-8pt]\nonumber
&&\hspace*{30pt}{} +\int^T_0\!\!\int
_E \bigl\llvert U_t(x) \bigr\rrvert
^2\nu^{\mathbb P}_t(dx)\,dt + \biggl(\int
_0^T \bigl\llvert\widehat F^{\mathbb P,0}_t
\bigr\rrvert \,dt \biggr)^2 \biggr],
\end{eqnarray}
for some constant $C_0$, independent of $\epsilon$. Now set $\epsilon
:=(2(1+C_0))^{-1}$ and plug (\ref{eqkk}) into (\ref{grahou}). One
then gets
\begin{eqnarray*}
&& \mathbb E^\mathbb P \biggl[\int_0^T
\bigl\llvert\hat a_t^{1/2}Z_t \bigr\rrvert
^2\,dt+\int_0^T\!\!\int
_EU^2_t(x)\nu^{\mathbb
P}_t(dx)\,dt
\biggr]
\\
&&\qquad \leq C\mathbb E^\mathbb P \biggl[\llvert\xi\rrvert
^2+ \sup_{0\leq t\leq T} \llvert Y_t\rrvert
^2 + \biggl(\int_0^T \bigl\llvert
\widehat F^{\mathbb P,0}_t \bigr\rrvert \,dt \biggr)^2
\biggr].
\end{eqnarray*}

From this and the estimate for $Y$, we immediately obtain
\[
\llVert Z\rrVert_{\mathbb H^{2,\kappa}_H}+\llVert U\rrVert_{\mathbb
J^{2,\kappa}_H}\leq C \bigl(
\llVert\xi\rrVert^2_{\mathbb L^{2,\kappa}_H}+\phi^{2,\kappa
}_H
\bigr).
\]

Then the estimate for $K^\P$ follows from (\ref{eqkk}). The estimates
for $z^\mathbb P$ and $u^\mathbb P$ can be proved similarly.
\end{pf}

%
\begin{Theorem}\label{estimates2}
Let Assumptions \ref{assumphref} and \ref{assumph2ref} hold. For
$i=1,2$, let us consider the solutions $(Y^i,Z^i,U^i,K^{\P,i})$ of the
2BSDEJ (\ref{2bsdej}) with terminal condition $\xi^i$. Then, there
exists a constant $C_\kappa$ depending only on $\kappa$, $T$ and the
Lipschitz constant of $F$ such that
\begin{eqnarray*}
&&  \bigl\llVert Y^1-Y^2 \bigr\rrVert_{\mathbb D^{2,\kappa}_H}
\leq C_{\kappa} \bigl\llVert\xi^1-\xi^2 \bigr \rrVert_{\mathbb L^{2,\kappa
}_H},
\\
&& \bigl\llVert Z^1-Z^2 \bigr\rrVert
^2_{\mathbb H^{2,\kappa}_H}+ \sup_{\mathbb P\in\mathcal
P^\kappa_H} \mathbb
E^\mathbb P \Bigl[\sup_{0\leq t\leq T} \bigl\llvert
K_t^{\P,1}-K_t^{\P,2} \bigr\rrvert
^2 \Bigr]+ \bigl\llVert U^1-U^2 \bigr\rrVert
^2_{\mathbb
J^{2,\kappa}_H}
\\
&&\qquad \leq C_{\kappa} \bigl\llVert\xi^1-\xi^2 \bigr
\rrVert_{\mathbb
L^{2,\kappa}_H} \bigl( \bigl\llVert\xi^1 \bigr\rrVert
_{\mathbb
L^{2,\kappa}_H}+ \bigl\llVert\xi^2 \bigr\rrVert_{\mathbb
L^{2,\kappa}_H}+
\bigl(\phi^{2,\kappa}_H \bigr)^{1/2} \bigr).
\end{eqnarray*}
%
\end{Theorem}

\begin{pf}
As in the previous theorem, we can obtain that there exists a constant
$C_\kappa$ depending only on $\kappa$, $T$ and the Lipschitz constant
of $\widehat F$, such that for all $\mathbb P$,
%
\begin{equation}
\label{estim2} \bigl\llvert y_t^{\mathbb P,1}-y_t^{\mathbb P,2}
\bigr\rrvert\leq C_\kappa\mathbb E_t^\mathbb P
\bigl[ \bigl\llvert\xi^1-\xi^2 \bigr\rrvert
^\kappa\bigr]^{1/\kappa}.
\end{equation}

Now by definition, we get from (\ref{estim2}) and representation
formula (\ref{representationref}) that
%
\begin{equation}
\label{eqy2} \bigl\llVert Y^1-Y^2 \bigr\rrVert
_{\mathbb D^{2,\kappa}_H}^2\leq C_\kappa\bigl\llVert
\xi^1-\xi^2 \bigr\rrVert^2_{\mathbb L^{2,\kappa}_H}.
\end{equation}

Applying It\^o's formula to $\llvert Y^1-Y^2\rrvert^2$, under
each $\mathbb
P\in\mathcal P^\kappa_H$, leads to
\begin{eqnarray*}
\nonumber
&& \mathbb E^\mathbb P \biggl[\int_0^T
\bigl\llvert\hat a_t^{1/2} \bigl(Z_t^1-Z_t^2
\bigr) \bigr\rrvert^2\,dt + \int_0^T\!\!\int_E \bigl\llvert U^1_t(x)-U^2_t(x)
\bigr\rrvert^2\nu^{\mathbb P}_t(dx)\,dt \biggr]
\\
&&\qquad \leq C\mathbb E^\mathbb P \bigl[ \bigl\llvert\xi^1-
\xi^2 \bigr\rrvert^2 \bigr]+\mathbb E^\mathbb P
\biggl[\int_0^T \bigl\llvert
Y_t^1-Y_t^2 \bigr\rrvert\, d
\bigl(K^{\mathbb
P,1}_t-K^{\mathbb P,2}_t \bigr)
\biggr]
\\
\nonumber
&&\quad\qquad{} +C\mathbb E^\mathbb P \biggl[\int_0^T
\bigl\llvert Y_t^1-Y_t^2 \bigr
\rrvert\biggl( \bigl\llvert Y_t^1-Y_t^2
\bigr\rrvert+ \bigl\llvert\hat a_t^{1/2}
\bigl(Z_t^1-Z_t^2 \bigr) \bigr
\rrvert
\\
&&\hspace*{140pt}{}+ \biggl(\int_E \bigl\llvert U^1_t(x)-U^2_t(x)
\bigr\rrvert^2\nu^{\mathbb P}_t(dx)\,dt
\biggr)^{1/2} \biggr)\,dt \biggr]
\\
&&\qquad \leq C \bigl( \bigl\llVert
\xi^1-\xi^2 \bigr\rrVert_{\mathbb L^{2,\kappa
}_H}^2+
\bigl\llVert Y^1-Y^2 \bigr\rrVert^2_{\mathbb D^{2,\kappa
}_H}
\bigr)
\\
\nonumber
&&\quad\qquad{} +\frac{1}2\mathbb E^\mathbb P \biggl[ \int
_0^T \bigl\llvert\hat a_t^{1/2}
\bigl(Z_t^1-Z_t^2 \bigr) \bigr
\rrvert^2\,dt+\int_0^T\!\!\int
_E \bigl\llvert U^1_t(x)-U^2_t(x)
\bigr\rrvert^2\nu^{\mathbb
P}_t(dx)\,dt \biggr]
\\
&&\quad\qquad{} +C \bigl\llVert Y^1-Y^2 \bigr\rrVert_{\mathbb D^{2,\kappa
}_H}
\Biggl(\mathbb E^\mathbb P \Biggl[\sum_{i=1}^2
\bigl(K^{\mathbb P,i}_T \bigr)^2 \Biggr]
\Biggr)^{1/2}.
\end{eqnarray*}

The estimates for $(Z^1-Z^2)$ and $(U^1-U^2)$ are now obvious from the
above inequality and the estimates of Theorem \ref{estimatesref}.
Finally the estimate for the difference of the nondecreasing processes
is obvious by definition.
\end{pf}

%
\begin{appendix}\label{app}
\section*{Appendix}
\setcounter{equation}{0}

\subsection{Generating and separable class of coefficients}\label{secdef}

We introduce the following notions inspired by \cite{stz3}:

\begin{Definitiona}
$\Ac_0 \subset\Ac_W$ is a generating class of coefficients if $\Ac
_0$ is stable for the concatenation operation; that is, if $(a,\nu),
(b,\beta) \in\Ac_0 \times\Ac_0$, then for each $t$,
\[
(a \mathbf{1}_{[0,t]} +b \mathbf{1}_{[t,+\infty)}, \nu\mathbf
{1}_{[0,t]} + \beta\mathbf{1}_{[t,+\infty)} ) \in\Ac_0.
\]
\end{Definitiona}

Notice that unlike \cite{stz3}, we do not impose their so-called
``constant disagreement time property,'' as it is only useful for them
to obtain their aggregation result, which, as mentioned before, is an
hopeless goal in our framework.

%
\begin{Definitiona}\label{sep}
We say that $\Ac$ is a separable class of coefficients generated by
$\Ac_0$ if $\Ac_0$ is a generating class of coefficients and if $\Ac
$ consists of all processes $a$ and random measures $\nu$ of the form
%
\begin{eqnarray}\label{doublesum}
a_t(\omega) &=& \sum_{n=0}^{+\infty}
\sum_{i=1}^{+\infty} a^{n,i}_t(
\omega) \mathbf{1}_{E^i_n}(\omega) \mathbf{1}_{[\tau
_n(\omega),\tau_{n+1}(\omega))}(t),
\nonumber\\[-8pt]\\[-8pt]\nonumber
\nu_t(\omega) &=& \sum_{n=0}^{+\infty}
\sum_{i=1}^{+\infty} \nu^{n,i}
_t(\omega)\mathbf{1}_{\Et^i_n}(\omega) \mathbf{1}_{[\tilde
{\tau}_n(\omega),\tilde{\tau}_{n+1}(\omega))}(t),
\end{eqnarray}
where for each $i$ and for each $n$, $(a^{n,i},\nu^{n,i}) \subset\Ac
_0$, $\tau_n$ and $\tilde{\tau}_n$ are $\F$-stopping times with
$\tau_0=0$, such that:
\begin{longlist}[(iii)]
\item[(i)] $\tau_n < \tau_{n+1}$ on $\{\tau_n < +\infty\}$
and $\widetilde\tau_n < \widetilde\tau_{n+1}$ on $\{\widetilde\tau
_n <
+\infty\}$.
\item[(ii)] $\inf\{n \geq0, \tau_n = +\infty\}+\inf\{n
\geq0, \widetilde\tau_n = +\infty\} <\infty$.
\item[(iii)] $\tau_n$ and $\widetilde\tau_n$ take countably
many values in some fixed $I_0\subset[0,T]$ which is countable and
dense in $[0,T]$.
\item[(iii)] For each $n$, $(E_i^n)_{i \geq1} \subset\Fc
_{\tau_n}$ and $(\widetilde{E}_i^n)_{i \geq1} \subset\Fc_{\tilde
{\tau
}_n}$ form a partition of $\Omega$.
\end{longlist}
\end{Definitiona}

%
\begin{Remarka}
If we refine the subdivisions, we can always take a common sequence of
stopping times $(\tau_n)_{n\geq0}$ and common sets $(E_i^n)_{i \geq
1,n\geq0}$ for $a$ and for $\nu$. This will be used throughout this section.
\end{Remarka}

The form for $a$ and $\nu$ in Definition \ref{sep} is directly
inspired by the so-called property of stability by concatenation and by
bifurcation in the theory of stochastic control. As shown, for instance,
in \cite{elk,elkarjean} or \cite{bt} (see Remark~3.1),
this property of control processes is tailor-made to be able to
retrieve the dynamic programming principle, and is somehow the minimal
stability property that must be verified. In our case, 2BSDEJs can be
seen formally as a weak version of a stochastic control problem for
which the controls are $a$ and $\nu$, and we will see below that the
set of probably measures we will consider will have this stability
property. This will be important for us in Proposition~4.2 of our
accompanying paper \cite{kpz}, where we recover the dynamic
programming principle.

The following proposition generalizes Proposition 4.11 of \cite
{stz3} and shows that a separable class of coefficients inherits the
``good'' properties of its generating class.

\begin{Propositiona}\label{toutrestevrai}
Let $\Ac$ be a separable class of coefficients generated by~$\Ac_0$. Then:
\begin{longlist}[(iii)]
\item[(i)] If $\Ac_0 \subset\Ac_W$, then $\Ac\subset\Ac_W$.
\item[(ii)] $\Ac$-quasi-surely is equivalent to $\Ac
_0$-quasi-surely, where for any $\tilde\Ac\subset\Ac_W$,
$\widetilde
\Ac$-q.s. means $\P$-a.s. for every $\P\in\{\P^\alpha_\nu,
(\alpha,\nu)\in\tilde\Ac\}$.
\item[(iii)] If every $\P\in\{ \P^{\alpha}_{\nu},
(\alpha, \nu) \in\Ac_0 \}$ satisfies the martingale
representation property, then every $\P\in\{ \P^{\alpha}_{\nu
}, (\alpha, \nu) \in\Ac\}$ also satisfies the martingale
representation property.
\item[(iv)] If every $\P\in\{ \P^{\alpha}_{\nu},
(\alpha, \nu) \in\Ac_0 \}$ satisfies the Blumenthal 0--1
law, then every probability measure $\P\in\{ \P^{\alpha}_{\nu
}, (\alpha, \nu) \in\Ac\}$ also satisfies the Blumenthal
0--1 law.
\end{longlist}
\end{Propositiona}

As in \cite{stz3}, to prove this result, we need the following two
lemmas. The first one is a straightforward generalization of Lemma~4.12
in \cite{stz3}, so we omit the proof. The second one is
analogous to Lemma~4.13 in \cite{stz3}.

%
\begin{Lemmaa} \label{simplesum}
Let $\Ac$ be a separable class of coefficients generated by $\Ac_0$.
For any $(a,\nu) \in\Ac$, and any $\F$-stopping time $\tau\in\Tc
$, there exist $\tilde{\tau} \in\Tc$ with $\tilde{\tau} \geq\tau
$, a~sequence $(a^n,\nu^n)_{n\geq1} \subset\Ac_0$ and a partition
$(E_n)_{n\geq1} \subset\Fc_{\tau}$ of $\Omega$ such that $\tilde
{\tau} > \tau$ on $\{\tau<+\infty\}$ and
%
\begin{eqnarray}\label{simplesumform}
\quad\qquad a_t(\omega) &=& \sum_{n\geq1}
a^n_t(\omega)\mathbf{1}_{E_n}(\omega)\quad
\mbox{and}\quad\nu_t(\omega) = \sum_{n\geq1}
\nu^n_t(\omega)\mathbf{1}_{E_n}(\omega),\qquad t<
\tilde{\tau}.
\end{eqnarray}

Finally, if $a$ and $\nu$ take the form (\ref{doublesum}) and $\tau
\geq\tau_n$, then we can choose $\tilde{\tau} \geq\tau_{n+1}$.
\end{Lemmaa}

\begin{pf}
We refer to the proof of Lemma 4.12 in \cite{stz3}.
\end{pf}

%
\begin{Lemmaa}\label{sumproba}
Let $\tau_1, \tau_2 \in\Tc$ be two stopping times such that $\tau
_1 \leq\tau_2$, and $(a^i,\nu^i)_{i\geq1} \subset\overline{\Ac
}_W$ and let $\{E_i, i\geq1\} \subset\Fc_{\tau_1}$ be a
partition of $\Omega$. Finally let $\P^0$ be a probability measure on
$\Fc_{\tau_1}$, and let $\{ \P^i, i\geq1 \}$ be a sequence of
probability measures such that for each $i$, $\P^i$ is a solution of
the martingale problem $(\P^0, \tau_1,\tau_2, a^i, \nu^i)$. Define
\begin{eqnarray*}
\P(E) &:=& \sum_{i \geq1} \P^i(E\cap
E_i)\qquad\mbox{for all } E\in\Fc_{\tau_2},
\\
a_t&:=& \sum_{i \geq1} a^i_t
\mathbf{1}_{E_i}\quad\mbox{and}\quad\nu_t:= \sum
_{i \geq1} \nu_t^i \mathbf{1}_{E_i},\qquad
t \in[\tau_1, \tau_2].
\end{eqnarray*}

Then $\P$ is a solution of the martingale problem $(\P^0, \tau
_1,\tau_2,a,\nu)$.
\end{Lemmaa}

\begin{pf}
By definition, $\P= \P^0$ on $\Fc_{\tau_1}$. In view of Remark \ref
{remmartpbm}, it is enough to prove that $M$, $J$ and $Q$ are $\P
$-local martingales on $[\tau_1, \tau_2]$. By localizing if
necessary, we may assume as usual that all these processes are actually
bounded. For any stopping times $\tau_1 \leq R \leq S \leq\tau_2$,
and any bounded $\Fc_R$-measurable random variable $\eta$, we have
\begin{eqnarray*}
\E^{\P} \bigl[ (M_S - M_R)\eta\bigr] &=&
\sum_{i\geq1} \E^{\P
^i} \bigl[(M_S - M_R)\eta\mathbf{1}_{E_i} \bigr]
\\
&=&  \sum
_{i\geq1} \E^{\P^i} \bigl[\E^{\P^i}
\bigl[(M_S-M_R)\mid\Fc_R \bigr]\eta
\mathbf{1}_{E_i} \bigr]
\\
&=& 0.
\end{eqnarray*}

Thus $M$ is a $\P$-local martingale on $[\tau_1, \tau_2]$. We can
prove similarly that $J$ and $Q$ are also $\P$-local martingales on
$[\tau_1, \tau_2]$.
\end{pf}

\begin{pf*}{Proof of Proposition \ref{toutrestevrai}}
The proof follows
closely the proof of Proposition 4.11 in \cite{stz3}, and we provide
it for the convenience of the reader.
\begin{longlist}[(i)]
\item[(i)] We take $(a,\nu) \in\Ac$. Let us prove that $(a,\nu)
\in\Ac_W$. We fix two stopping times $\theta_1,\theta_2$ in $\Tc$.
We define a sequence $(\tilde{\tau}_n)_{n\geq0}$ as follows:
\[
\tilde{\tau}_0:= \theta_1\quad\mbox{and}\quad\tilde{
\tau}_n:= (\tau_n \vee\theta_1) \wedge
\theta_2,\qquad n\geq1.
\]
To prove that the martingale problem $(\P^0, \theta_1, \theta
_2,a,\nu)$ has a unique solution, we prove by induction on $n$ that
the martingale problem $(\P^0, \tilde{\tau}_0,\tilde{\tau}_n,a,\nu
)$ has a unique solution.
\begin{longlist}
\item[\textit{Step} 1 \textit{of the induction}.]
Let $n=1$, and let us first construct a
solution to the martingale problem $(\P^0, \tilde{\tau}_0, \tilde
{\tau}_1,a,\nu)$. For this purpose, we apply Lemma \ref{simplesum}
with $\tau= \tilde{\tau}_0$ and $\tilde{\tau}=\tilde{\tau}_1$,
which leads to $a_t = \sum_{i\geq1} a_t^i \mathbf{1}_{E_i}$ and $\nu
_t = \sum_{i\geq1} \nu^i_t \mathbf{1}_{E_i}$ for all $t<\tilde
{\tau}_1$, where $(a^i,\nu^i) \in\Ac_0$ and $\{E_i, i\geq1\}
\subset\Fc_{\tilde{\tau}_0}$ forms a partition of $\Omega$. For $i
\geq1$, let $\P^{0,i}$ be the unique solution of the martingale
problem $(\P^0, \tilde{\tau}_0, \tilde{\tau}_1,a_i,\nu_i)$ and define
\[
\P^{0,a}(E):= \sum_{i \geq1}
\P^{0,i}(E \cap E_i)\qquad\mbox{for all } E \in
\Fc_{\tilde{\tau}_1}.
\]

Then Lemma \ref{sumproba} tells us that $\P^{0,a}$ solves the
martingale problem $(\P^0, \tilde{\tau}_0,\break  \tilde{\tau}_1, a,\nu
)$. Now let $\P$ be an arbitrary solution of the martingale problem
$(\P^0, \tilde{\tau}_0, \tilde{\tau}_1, a,\nu)$, and let us prove
that $\P= \P^{0,a}$. We first define
\[
\P^i(E):= \P(E\cap E_i) + \P^{0,i} \bigl(E
\cap E_i^c \bigr)\qquad \forall E \in\Fc_{\tilde{\tau}_1}.
\]

Using Lemma \ref{sumproba}, and the facts that $a^i =a \mathbf
{1}_{E_i} + a^i \mathbf{1}_{E_i^c}$ and $\nu^i =\nu\mathbf{1}_{E_i}
+ \nu^i \mathbf{1}_{E_i^c}$, we conclude that $\P^i$ solves the
martingale problem $(\P^0, \tilde{\tau}_0, \tilde{\tau}_1,a^i,\nu
^i)$. Since\vspace*{1pt} this problem has a unique solution, we thus have $\P^i =
\P^{0,i}$ on $\Fc_{\tilde{\tau}_1}$. This implies that for each $i
\geq1$ and for each $E \in\Fc_{\tilde{\tau}_1}$, $\P^i(E\cap E_i)
= \P^{0,i}(E \cap E_i)$, and finally
\[
\P^{0,a}(E) = \sum_{i \geq1}
\P^{0,i}(E \cap E_i) = \sum_{i \geq
1}
\P^{i}(E \cap E_i) = \P(E)\qquad \forall E \in
\Fc_{\tilde{\tau}_1}.
\]
\end{longlist}
\begin{longlist}
\item[\textit{Step} 2 \textit{of the induction}.] We assume that the 
martingale problem
$(\P^0, \tilde{\tau}_0,\break \tilde{\tau}_n, a,\nu)$ has a unique
solution denoted by $\P^n$. Using the same reasoning as above, we see
that the martingale problem $(\P^n, \tilde{\tau}_n,\tilde{\tau
}_{n+1},a,\nu)$ has a unique solution, denoted by $\P^{n+1}$. Then
the processes $M$, $J$ and $Q$ defined\vspace*{1pt} in Remark \ref{remmartpbm}
are $\P^{n+1}$-local martingales on $[\tilde{\tau}_n,\tilde{\tau
}_{n+1}]$, and since $\P^{n+1}$ coincides with $\P^{n}$ on $\Fc
_{\tilde{\tau}_n}$, $M$, $J$ and $Q$ are also $\P^{n+1}$-local
martingales on $[\tilde{\tau}_0,\tilde{\tau}_n]$. Hence $\P^{n+1}$
solves the martingale problem $(\P^0, \tilde{\tau}_0,\tilde{\tau
}_{n+1},a,\nu)$. We suppose now that $\P$ is another arbitrary
solution to the martingale problem $(\P^0, \tilde{\tau}_0,\tilde
{\tau}_{n+1},a,\nu)$. By the induction assumption, $\P^n = \P$ on
$\Fc_{\tilde{\tau}_n}$. Then $\P$ solves the martingale problem
$(\P^n, \tilde{\tau}_n,\tilde{\tau}_{n+1},a,\nu)$, and by
uniqueness $\P= \P^{n+1}$ on $\Fc_{\tilde{\tau}_{n+1}}$. The
induction is now complete.

Note that $\Fc_{\theta_2} = \bigvee_{n \geq1} \Fc_{\tilde{\tau
}_n}$. Indeed, since $\inf\{ n\geq1, \tau_n = +\infty\} < +\infty
$, then $\inf\{ n\geq1, \tilde{\tau}_n = \theta_2 \} < +\infty$.
This allows us to define $\P^{\infty}(E):= \P^n(E)$ for $E \in\Fc
_{\tilde{\tau}_n}$ and to extend it uniquely to $\Fc_{\theta_2}$.
Now using again Remark \ref{remmartpbm}, we conclude that $\P
^{\infty}$ solves $(\P^0, \theta_1, \theta_2, a, \nu)$ and is unique.
\end{longlist}

\item[(ii)] We now prove that $\Ac$-quasi-surely is equivalent to
$\Ac_0$-quasi-surely.

We take $(a,\nu) \in\Ac$ and we apply Lemma \ref{simplesum} with
$\tau= + \infty$ to write $a_t = \sum_{i\geq1} a^i_t \mathbf
{1}_{E_i}$ and $\nu_t = \sum_{i\geq1} \nu^i_t \mathbf{1}_{E_i}$
for all $t\geq0$, where $(a^i,\nu^i) \in\Ac_0$ and $\{E_i, i\geq
1\} \subset\Fc_{\infty}$ forms a partition of $\Omega$. Take a set
$E$ such that $\P^{\tilde{a}}_{\tilde{\nu}}(E)=0$ for every
$(\tilde{a},\tilde{\nu}) \in\Ac_0$, then
\[
\P^a_{\nu}(E) = \sum_{i \geq1}
\P^a_{\nu}(E \cap E_i) = \sum
_{i
\geq1} \P^{a^i}_{\nu^i}(E \cap
E_i) = 0.
\]

\item[(iii)] By (i), since $\Ac\subset\Ac_W$, for any $(a,\nu)\in
\Ac$, the corresponding martingale problem has a unique solution,
which is therefore an extremal point in the set of solutions. Hence, we
can apply Theorem III.4.29 in \cite{jac} to obtain immediately the
predictable representation property.

\item[(iv)] Take $(a,\nu) \in\Ac$ of the form (\ref{doublesum}),
in which we can take $\tau_0=0$ without loss of generality. $\P
^a_{\nu}$ is the law on $[0,\tau_1]$ of a semimartingale with
characteristics
\[
\biggl( -\int_{0}^{t} \int_E
x \mathbf{1}_{\llvert x\rrvert
>1}\tilde{\nu}_s(dx)\,ds, \int
_{0}^{t} \tilde{a}_s \,ds, \tilde{\nu
}_s(dx)\,ds \biggr),
\]
where
\[
\tilde{a}_t:= \sum_{i \geq1}
a^{0,i} \mathbf{1}_{E_0^i}\quad\mbox{and}\quad\tilde{
\nu}_t:= \sum_{i \geq1} \nu^{0,i}
\mathbf{1}_{E_0^i},
\]
where $\{E_0^i, i \geq1\} \subset\Fc_0$ is a partition of $\Omega
$. Since $\Fc_0$ is trivial, the partition is only composed of $\Omega
$ and $\varnothing$, and then
\[
\tilde{a}_t:= a_t^{0,1}\quad\mbox{and}\quad
\tilde{\nu}_t = \nu^{0,1}_t.
\]

Then for $E \in\Fc_{0^+}$, $\P^a_{\nu}(E) = \P^{\tilde
{a}}_{\tilde{\nu}}(E) = 0$,
because $\P^{\tilde{a}}_{\tilde{\nu}}$ satisfies the Blumenthal
0--1 law by hypothesis.\quad\qed
\end{longlist}\noqed
\end{pf*}

\subsection{The measures \texorpdfstring{$\mathbb{P}^{\alpha,\beta}_F$}{Palpha,betaF}}

\begin{Lemmaa}
\label{upward}
Fix an arbitrary measure $\P=\P^{\alpha,\beta}_F$ in $\mathcal
{P}_H^{\kappa}$. The set $\mathcal{P}_H^{\kappa}(t^+,\P)$ is upward
directed; that is, for each $\P_1:= \P^{\alpha_1,\beta_1}_{F_1}$
and $\P_2:= \P^{\alpha_2,\beta_2}_{F_2}$ in $\mathcal{P}_H^{\kappa
}(t^+,\P)$, there exists $\P' \in\mathcal{P}_H^{\kappa}(t^+,\P)$
such that $\forall u>t$,
%
\begin{eqnarray}\label{defupdir}
\qquad && \E_t^{\P'} \bigl[ \bigl(K^{\P'}_u
- K^{\P'}_t \bigr)^2 \bigr] = \max\bigl\{
\E_t^{\P_1} \bigl[ \bigl(K^{\P_1}_u -
K^{\P _1}_t \bigr)^2 \bigr],
\E_t^{\P_2} \bigl[ \bigl(K^{\P_2}_u -
K^{\P_2}_t \bigr)^2 \bigr] \bigr\}.
\end{eqnarray}
\end{Lemmaa}

\begin{pf}
We define the following $\mathcal{F}_{t^+}$-measurable sets:
\[
E_1:= \bigl\{ \omega\in\Omega\dvtx \E_t^{\P_2}
\bigl[ \bigl(K^{\P
_2}_u - K^{\P_2}_t
\bigr)^2 \bigr](\omega) \leq\E_t^{\P
_1} \bigl[
\bigl(K^{\P_1}_u - K^{\P_1}_t
\bigr)^2 \bigr](\omega) \bigr\},
\]
and $E_2:= \Omega\setminus E_1$. Then for all $A \in
\mathcal{F}_T$, we define the probability measure $\P'$ by
\[
\P'(A):= \P_1(A\cap E_1) +
\P_2(A \cap E_2).
\]
By definition, $\P'$
satisfies (\ref{defupdir}). Let us prove now that $\P' \in\mathcal
{P}_H^{\kappa}(t^+,\P)$. As in the proof of claim (4.17) in \cite
{stz}, for $s \in[0,T]$, we define the processes $\alpha^*$, $\beta
^*$ and the measure $F^*$ as follows:
\begin{eqnarray*}
\alpha^*_s(\omega)&:=& \alpha_s(\omega)
\mathbf{1}_{[0,t)}(s)
\\
&&{} + \bigl( \alpha_s^1(
\omega)\mathbf{1}_{ \{X^{\alpha,\beta} \in
E_1 \}}(\omega) + \alpha_s^2(
\omega)\mathbf{1}_{ \{
X^{\alpha,\beta} \in E_2 \}}(\omega) \bigr)\mathbf{1}_{[t,T]}(s),
\\
\beta^*_s(\omega,x)&:=& \beta_s(\omega,x)
\mathbf{1}_{[0,t)}(s)
\\
&&{} + \bigl( \beta_s^1(
\omega,x)\mathbf{1}_{ \{X^{\alpha,\beta}
\in E_1 \}}(\omega) + \beta_s^2(
\omega,x)\mathbf{1}_{ \{
X^{\alpha,\beta} \in E_2 \}}(\omega) \bigr)\mathbf{1}_{[t,T]}(s),
\\
F^*_s(\omega)&:=& F_s(\omega) \mathbf{1}_{[0,t)}(s)
\\
&&{}
+ \bigl( F_s^1(\omega)\mathbf{1}_{ \{X^{\alpha,\beta} \in E_1 \}
}(
\omega) + F_s^2(\omega)\mathbf{1}_{ \{X^{\alpha,\beta} \in
E_2 \}}(
\omega) \bigr)\mathbf{1}_{[t,T]}(s),
\end{eqnarray*}
where $X^{\alpha,\beta}$ is defined in (\ref{Xalphabeta}).

First of all, we clearly have $F^*\in\Vc$, since this set is stable
by concatenation and bifurcation by definition. We can therefore define
the probability measure $\mathbb P_{0,F^*}$.\footnote{The attentive
reader may have remarked that $F^*$ is not defined for every $\omega$,
but only for those such that their path up to time $t^+$ is in the
support of $\P$ restricted to $\Fc_{t^+}$. This may appear as a
problem, however, since we know that the measure $\P^{\alpha^*, \beta
^*}$ has to agree with $\P$ on $\Fc_{t^+}$, we actually only need to
solve the martingale problem in the definition of $\P_{0,F^*}$
starting from time $t$.} Moreover, we have
\[
0< \underline{\alpha} \wedge\underline{\alpha}^1 \wedge\underline{
\alpha}^2 \leq\alpha^* \leq\overline{\alpha} \vee\overline{
\alpha}^1 \vee\overline{\alpha}^2,
\]
where $\underline
{\alpha}$, $\overline{\alpha}$, $\underline{\alpha}^i$, $\overline
{\alpha}^i$ are the lower and upper bounds of the processes $\alpha$
and $\alpha^i$. Next, we have to check that $\beta^*\in\Rc_{F^*}$.
It is clear that for every $\omega\in\Omega$, $F^*(dx)$-a.e.,
\begin{eqnarray*}
\bigl\llvert\beta^*_s \bigr\rrvert(\omega,x)&\leq&\bigl(C{\mathbf
1}_{0\leq s< t}+ \bigl(C_1{\mathbf1}_{X^{\alpha,\beta}\in E_1}(
\omega)+C_2{\mathbf1}_{X^{\alpha,\beta}\in E_2}(\omega) \bigr){\mathbf
1}_{t\leq s\leq T}
\bigr) \bigl(1\wedge\llvert x\rrvert\bigr)
\\
&\leq& C^* \bigl(1\wedge\llvert x\rrvert\bigr),
\end{eqnarray*}
since $F^*$ coincides with $F$ before $t$ and with either $F^1$ or
$F^2$ after $t$.

The strict monotony of $x\longmapsto\beta_s^*(\omega,x)$ for
Lebesgue almost every $s\in[0,T]$ and $\P_{0,F^*}$-a.e. $\omega\in
\Omega$ follows similarly from the corresponding properties of $\beta
$, $\beta^1$ and $\beta^2$ and the fact that the support of the law
of the jumps of $B$ at time $s$ under $\P_{0,F^*}$ coincides with the
support of the same law under $\P_{0,F}$ for $s<t$ and under either
$\P_{0,F^1}$ or $\P_{0,F^2}$ for $s\geq t$.

We can check similarly that
\[
\int_0^T\!\!\int_{\{\llvert x\rrvert>1\}} x
\nu^{F^*,\beta
^*}_s(dx,ds)<+\infty
\]
and
\[
\mathbb
E^\P\biggl[\int_0^T\!\!\int
_E\llvert x\rrvert^2\nu^{F^*,\beta^*}_s(dx)\,ds
\biggr]<+\infty.
\]

Therefore, we have proved that $\P^{\alpha^*,\beta^*}_{F^*}\in
\overline{\Pc}_S$. Moreover, using the same arguments as in the step~3
of the proof of Lemma~A.3 in \cite{kpz}, we can easily show that
$\P'=\P^{\alpha^*, \beta^*}_{F^*}$. Finally, we compute
\begin{eqnarray*}
&& \E^{\P'} \biggl[ \int_0^T \bigl
\llvert\widehat{F}^{\P',0}_s \bigr\rrvert^2 \,ds
\biggr]
\\
&&\qquad =\E^{\P} \biggl[ \int_0^t
\bigl\llvert\widehat{F}^{\P,0}_s \bigr\rrvert
^2 \,ds \biggr] + \E^{\P_1} \biggl[ \int_t^T
\bigl\llvert\widehat{F}^{\P
_1,0}_s \bigr\rrvert
^2 \,ds\, \mathbf{1}_{E_1} \biggr]
\\
&&\quad\qquad{} +\E^{\P_2} \biggl[
\int_t^T \bigl\llvert\widehat
{F}^{\P_2,0}_s \bigr\rrvert^2 \,ds\,
\mathbf{1}_{E_2} \biggr]
\\
&&\qquad \leq\E^{\P} \biggl[ \int_0^T \bigl
\llvert\widehat{F}^{\P,0}_s \bigr\rrvert^2 \,ds
\biggr] + \E^{\P_1} \biggl[ \int_0^T
\bigl\llvert\widehat{F}^{\P_1,0}_s \bigr\rrvert
^2 \,ds\, \mathbf{1}_{E_1} \biggr]
\\
&&\quad\qquad{}+\E^{\P_2} \biggl[
\int_0^T \bigl\llvert\widehat
{F}^{\P_2,0}_s \bigr\rrvert^2 \,ds\,
\mathbf{1}_{E_2} \biggr].
\\
&&\qquad < +\infty.
\end{eqnarray*}
Since by construction $\P'$ coincides with $\P$ on $\mathcal
{F}_{t^+}$, we have indeed shown that $\P'\in\Pc^\kappa_H(t^+,\P)$.
\end{pf}

\subsection{$L^r$-integrability of exponential martingales}

%
\begin{Lemmaa}\label{inegfonc1}
Let $\delta>0$ and $n\in\mathbb N^*$. Then there exists a constant
$C_{n,\delta}$ depending only on $\delta$ and $n$ such that
\[
(1+x)^{-n}-1+nx\leq C_{n,\delta}x^2\qquad\mbox{for
all }x\in[-1+\delta,+\infty).
\]
\end{Lemmaa}

\begin{pf}
The inequality is clear for $x$ large enough; let us say $x\geq M$ for
some $M>0$. Then a simple Taylor expansion shows that this also holds
in a neighborhood of $0$, that is to say for $x\in[-\varepsilon
,\varepsilon]$ for
some $\varepsilon>0$. Finally, for $x\in(-1+\delta,-\varepsilon
)\cup(\varepsilon,M)$,
it is clear that we can choose $C$ large enough such that the
inequality also holds.
\end{pf}

M\'emin \cite{memin} and then L\'epingle and M\'emin \cite
{lepinmemin2} proved some useful multiplicative decompositions of
exponential semimartingales. We give here one of these representations
that we will use in the proof of Lemma \ref{expomart1}.


\begin{Propositiona}[(Proposition II.1 of \cite{lepinmemin1})]\label{proplepinmemin}
Let $N$ be a local martingale and let $A$ be a predictable process with
finite variation such that $\Delta A \neq-1$. We assume $N_0=A_0=0$.
Then there exists a local martingale $\widetilde{N}$ with $\widetilde{N}_0=0$
and such that
\[
\mathcal{E}(N+A) = \mathcal{E}(\widetilde{N})\mathcal{E}(A).
\]
\end{Propositiona}

%
\begin{Lemmaa}\label{expomart1}
Let $\lambda>0$ and $M$ be a local martingale with bounded jumps, such
that $\Delta M\geq-1+\delta$, for a fixed $\delta>0$. Let
$V^{-\lambda}$ be the predictable compensator of
\[
\biggl\{ W_t^{-\lambda} = \sum_{s\leq t}
\bigl[(1+\Delta M_s)^{-\lambda
}-1+\lambda\Delta
M_s \bigr], t \geq0 \biggr\}.
\]
We have:
\begin{longlist}[(ii)]
\item$\mathcal{E}^{-\lambda}(M) = \mathcal{E}(N^{-\lambda} +
A^{-\lambda})$ where
\[
A^{-\lambda} = \frac{\lambda(\lambda
+1)}{2} \bigl\langle M^c,M^c
\bigr\rangle^T + V^{-\lambda},\qquad N^{-\lambda} = -\lambda
M^T + W^{-\lambda} - V^{-\lambda}.
\]
\item There exist a local martingale $\widetilde{N}^{-\lambda}$ such that
\[
\mathcal{E}^{-\lambda}(M) = \mathcal{E} \bigl(\widetilde{N}^{-\lambda}
\bigr) \mathcal{E} \bigl(A^{-\lambda} \bigr).
\]
\end{longlist}
\end{Lemmaa}

\begin{pf}
First note that thanks to Lemma \ref{inegfonc1}, for $\lambda>0$,
$(1+x)^{-\lambda}-1+\lambda x \leq C x^2$, and thus $W^{-\lambda}$ is
integrable. We set
\[
T_n = \inf\biggl\{t\geq0\dvtx \mathcal{E}(M)_t \leq
\frac{1}{n} \biggr\} \quad\mbox{and}\quad M^n_t =
M_{t \wedge T_n}.
\]
Then $M^n$ and $ \mathcal{E}(M^n)$ are local martingales, $\mathcal
{E}(M^n) \geq\frac{1}{n}$ and $\mathcal{E}(M^n)_t = \mathcal
{E}(M)_t$ if $t<T_n$. The assumption $\Delta M >-1$ shows that $T_n$
tends to infinity when $n$ tends to infinity. For each $n\geq1$, we
apply It\^o's formula to a $\mathcal{C}^2$ function $f_n$ that
coincides with $x^{-\lambda}$ on $ [\frac{1}{n}, +\infty)$,
\begin{eqnarray*}
\mathcal{E}^{-\lambda} \bigl(M^n \bigr)_t &=& 1 -
\lambda\int_0^t \mathcal{E}^{-\lambda-1}
\bigl(M^n \bigr)_{s^-} \,d \mathcal{E} \bigl(M^n
\bigr)_s
\\
&&{}+ \frac{\lambda
(\lambda+1)}{2} \int_0^t
\mathcal{E}^{-\lambda-2} \bigl(M^n \bigr)_{s^-} \,d \bigl
\langle\bigl(\mathcal{E} \bigl(M^n \bigr) \bigr)^c \bigr
\rangle_s
\\
&&{} + \sum_{s\leq t} \bigl[ \mathcal{E}^{-\lambda}
\bigl(M^n \bigr)_s - \mathcal{E}^{-\lambda}
\bigl(M^n \bigr)_{s^-} + \lambda\mathcal{E}^{-\lambda
-1}
\bigl(M^n \bigr)_{s^-} \Delta\mathcal{E}^{-\lambda}
\bigl(M^n \bigr)_s \bigr]
\\
&=& 1+ \int_0^t \mathcal{E}
\bigl(M^n \bigr)_{s^-} \,dX^n_s,
\end{eqnarray*}
where
\[
X^n_t:= -\lambda M^n_t +
\frac{\lambda(\lambda+1)}{2} \bigl\langle\bigl(M^n \bigr)^c,
\bigl(M^n \bigr)^c \bigr\rangle_t +\sum
_{s\leq t} \bigl[(1+\Delta M_s)^{-\lambda}-1+
\lambda\Delta M_s \bigr],
\]
and then $\mathcal{E}^{-\lambda}(M^n) = \mathcal{E}(X^n)$. Let us
define the nontruncated counterpart $X$ of $X^n$:
\[
X = -\lambda M + \frac{\lambda(\lambda+1)}{2} \bigl\langle M^c,
M^c \bigr\rangle+ W^{-\lambda}.
\]
On the interval $[0,T_n[$, we have $X^n=X$ and $\mathcal{E}^{-\lambda
}(M) = \mathcal{E}(X)$, now letting $n$ tend to infinity, we obtain
that $\mathcal{E}^{-\lambda}(M)$ and $\mathcal{E}(X)$ coincide on
$[0,+\infty[$, which is the point (i) of the lemma.

We want to use Proposition \ref{proplepinmemin} to prove the point
(ii), so we need to show that $\Delta A >-1$. We set
\[
S= \inf\bigl\{t \geq0\dvtx \Delta A_t^{-\lambda} \leq-1 \bigr\}.
\]
It is a predictable stopping time. Using this, and the fact that $M$
and $(W^{-\lambda} - V^{-\lambda})$ are local martingales, we have
\[
\Delta A_S^{-\lambda} = \E\bigl[\Delta A_S^{-\lambda}
\mid\mathcal{F}_{S^-} \bigr] =\E[\Delta X_S \mid
\mathcal{F}_{S^-} ]= \E\bigl[(1+\Delta M_S)^{-\lambda}
\mid\mathcal{F}_{S^-} \bigr],
\]
and since $\{S<+\infty\} \in\mathcal{F}_{S^-}$,
\[
0 \geq\E\bigl[ \mathbf{1}_{\{S<+\infty\}} \bigl(1+\Delta A_S^{-\lambda
}
\bigr) \bigr] = \E\bigl[ \mathbf{1}_{\{S<+\infty\}} (1+\Delta
M_S)^{-\lambda} \bigr].
\]
Then $\Delta M_S \leq-1$ on $\{S<+\infty\}$, which means that
$S=+\infty$ and $\Delta A >-1$ a.s. The proof is now complete.
\end{pf}

We are finally in a position to state the lemma on $L^r$ integrability
of exponential martingales for a negative exponent $r$.

%
\begin{Lemmaa}\label{expomart2}
Let $\lambda>0$ and let $M$ be a local martingale with bounded jumps,
such that $\Delta M\geq-1+\delta$, for a fixed $\delta>0$, and
$\langle M,M\rangle_t$ is bounded $dt \times\P$-a.s. Then
\[
\E^{\P} \bigl[ \mathcal{E}(M)_t^{-\lambda} \bigr] < +
\infty \qquad dt \times\P\mbox{-a.s.}
\]
\end{Lemmaa}

\begin{pf}
Let $n \geq1$ be an integer. We will denote $\tilde{\mu}_M = \mu_M
- \nu_M$ the compensated jump measure of $M$. Thanks to Lemma \ref
{expomart1}, we write the decomposition
\[
\mathcal{E}(M)^{-n} = \mathcal{E} \bigl(\widetilde{N}^{-n}
\bigr)\mathcal{E} \bigl(\tfrac{1}{2} n(n+1) \bigl\langle M^c,
M^c \bigr\rangle+V^{-n} \bigr),
\]
where $\widetilde{N}^{-n}$ is a local martingale, and $V^{-n}$ is
defined as
$V^{-\lambda}$.
Using Lemma \ref{inegfonc1}, we have the inequality
\[
V^{-n}_t \leq\int_0^t\!\int_E C x^2 \nu_M(dx,ds),
\]
and using the previous representation we obtain
\begin{eqnarray*}
\mathcal{E}(M)^{-n}_t &\leq&\mathcal{E} \bigl(
\widetilde{N}^{-n} \bigr)_t \mathcal{E} \biggl(
\frac{1}{2} n(n+1) \bigl\langle M^c, M^c \bigr
\rangle+\int_0^{\cdot} \int_E
C x^2 \nu_M(dx,ds) \biggr)_t
\\
&\leq& \mathcal{E} \bigl(
\widetilde{N}^{-n} \bigr)_t \exp\biggl( \biggl(
\frac{1}{2} n(n+1)+C \biggr) \langle M, M\rangle_t \biggr)
\\
&\leq& C \mathcal{E} \bigl(\widetilde{N}^{-n} \bigr)_t
\qquad\mbox{since } \langle M, M\rangle_t\mbox{ is bounded.}
\end{eqnarray*}

Let us prove now that the jumps of $\widetilde{N}^{-n}$ are strictly
bigger than $-1$. We compute
\begin{eqnarray*}
\Delta\widetilde{N}^{-n} &=& \frac{ \Delta N^{-n}}{1+ \Delta
A^{-n}}\qquad\mbox{where }
A^{-n}\mbox{ is defined as in Lemma \ref
{expomart2}}
\\
&=& \frac{(1+\Delta M)^{-n}}{1+\Delta V^{-n}} -1 >-1\qquad\mbox{since }
-1< \Delta M \leq B
\mbox{ and }\Delta V^{-n} >-1.
\end{eqnarray*}

This implies that $\mathcal{E}(\widetilde{N}^{-n})$ is a positive
supermartingale which equals $1$ at $t=0$. We deduce
\[
\E\bigl[ \mathcal{E}(M)^{-n}_t \bigr] \leq C \E\bigl[
\mathcal{E} \bigl(\widetilde{N}^{-n} \bigr)_t \bigr] \leq
C.
\]
We have the desired integrability for negative integers. We extend the
property to any negative real number by H\"older's inequality.
\end{pf}
\end{appendix}

\section*{Acknowledgments}
Part of this work was carried out while all the all authors
were working at CMAP, Ecole Polytechnique.
Part of this work was also carried out while Dylan Possama\"i was invited to
the Mathematics Department at the National
University of Singapore. The authors would like to thank warmly Marcel Nutz for his precious
advices and for pointing out
a mistake in a previous version of the paper. They are also grateful to
an Associate Editor and
two referees who greatly helped to improve the readability of the paper.


%

\printaddresses
\end{document}